\documentclass[11pt]{article}

\newcommand{\documentdate}{27 February, 2019}

\usepackage{a4wide,latexsym,graphicx,varioref,color,amsmath}

\vspace{-6.0in}
\title{
  High-Order Evaluation Complexity for Convexly-Constrained Optimization with Non-Lipschitzian Group Sparsity Terms
  }
\author{\small
  X. Chen\thanks{
    Department of Applied Mathematics, The Hong Kong Polytechnic University, Hong Kong.
    Email: \tt{xiaojun.chen@polyu.edu.hk}
  }
~and  Ph. L. Toint\thanks{
    Namur Centre for Complex Systems (naXys),
    University of Namur,
    61, rue de Bruxelles, B-5000 Namur, Belgium.
    Email: \tt{philippe.toint@unamur.be}}
}
\date{\small \documentdate}
\newcommand{\ass}[2]{\label{ass-#1}
                     \begin{list}{}{\setlength{\leftmargin}{2.5cm}}
                     \item \hspace{-2.6cm} \framebox[2.0cm]{\bf #1} \,\,#2
                     \end{list}}

\DeclareMathOperator*\spanset{span}

\newcommand{\comment}[1]{}
\newcommand{\numsection}[1]{\section{#1}\setcounter{equation}{0}}
\newcommand{\appnumsection}[1]{\section*{#1}\setcounter{equation}{0}
  \renewcommand{\theequation}{A.\arabic{equation}}
  \renewcommand{\thetheorem}{A.\arabic{theorem}}
  \renewcommand{\thetable}{A.\arabic{table}}
  \renewcommand{\thefigure}{A.\arabic{figure}}
  \renewcommand{\thesection}{A} }
\renewcommand{\theequation}{\arabic{section}.\arabic{equation}}
\renewcommand{\thetable}{\arabic{section}.\arabic{table}}
\renewcommand{\thefigure}{\arabic{section}.\arabic{figure}}

\newcommand{\calA}{{\cal A}}
\newcommand{\calC}{{\cal C}} \newcommand{\calR}{{\cal R}}
\newcommand{\calN}{{\cal N}} \newcommand{\calH}{{\cal H}}
\newcommand{\calX}{{\cal X}} \newcommand{\calM}{{\cal M}}
\newcommand{\calW}{{\cal W}} \newcommand{\calF}{{\cal F}}
\newcommand{\calI}{{\cal I}} 
 \newcommand{\calV}{{\cal V}}
\newcommand{\calO}{{\cal O}} \newcommand{\calS}{{\cal S}}
\newcommand{\calU}{{\cal U}} 
 \newcommand{\calL}{{\cal L}}
\newcommand{\eqdef}{\stackrel{\rm def}{=}}
\newcommand{\req}[1]{(\ref{#1})}
\newcommand{\beqn}[1]{\begin{equation}\label{#1}}
\newcommand{\eeqn}{\end{equation}}
\newcommand{\tim}[1]{\;\; \mbox{#1} \;\;}
\newcommand{\ms}{\;\;\;\;}

\newcommand{\mystack}[2]{_{\stackrel{\scriptstyle #1}{\scriptstyle #2}}}
\newcommand{\bpr}{{\bf Proof.} \hspace{0.2mm}}
\newcommand{\epr}{\hfill $\Box$ \vspace*{1em}}
\newcommand{\proof}[1]{
\begin{list}{}{
\setlength{\topsep}{0.0pt}
\setlength{\partopsep}{0.0pt}
\setlength{\leftmargin}{0.0\textwidth}
\setlength{\rightmargin}{0.0\leftmargin}
\setlength{\labelwidth}{0.0\leftmargin}
\setlength{\labelsep}{0.0\leftmargin}}
\item \bpr #1 \epr \noindent
\end{list}}
\renewcommand{\Re}{\hbox{I\hskip -2pt R}}

\newcommand{\Na}{\hbox{I\hskip -1.8pt N}}

\newcommand{\sfrac}[2]{{\scriptstyle \frac{#1}{#2}}}
\newcommand{\half}{\sfrac{1}{2}}

\newcommand{\quarter}{\sfrac{1}{4}}
\newcommand{\threequarters}{\sfrac{3}{4}}
\newcommand{\bigsum}{\displaystyle \sum}
\newcommand{\bigmax}{\displaystyle \max}
\newcommand{\bigmin}{\displaystyle \min}

\newcommand{\bigfrac}[2]{\frac{\displaystyle #1}{\displaystyle #2}}

\newcommand{\kap}[1]{\kappa_{\mbox{\rm \tiny #1}}}
\newcommand{\ii}[1]{\{1, \ldots, #1 \}}
\newcommand{\iibe}[2]{\{ #1, \ldots, #2 \}}

\newtheorem{theorem}{Theorem}[section]
\newtheorem{lemma}[theorem]{Lemma}
\newtheorem{corollary}[theorem]{Corollary}

\newcommand{\lthm}[2]{\begin{theorem} \label{#1} {\it #2} \end{theorem} }
\newcommand{\llem}[2]{\begin{lemma} \label{#1} {\it #2} \end{lemma} }

\newcounter{algo}[section]
\renewcommand{\thealgo}{\thesection.\arabic{algo}}

\newcommand{\algo}[3]{\refstepcounter{algo}
\begin{center}\begin{figure}[htbp]
\framebox[\textwidth]{
\parbox{0.95\textwidth} {\vspace{\topsep}
{\bf Algorithm \thealgo : #2}\label{#1}\\
\vspace*{-\topsep} \mbox{ }\\
{#3} \vspace{\topsep} }}
\end{figure}\end{center}}

\begin{document}

\maketitle

\vspace*{-5mm}
\begin{abstract}
{\small This paper studies high-order evaluation complexity for partially separable convexly-constrained optimization involving non-Lipschitzian group sparsity terms in a nonconvex objective function. We propose a partially separable adaptive regularization algorithm using a $p$-th order Taylor model and show that the algorithm can produce an
$(\epsilon,\delta)$-approximate $q$-th-order stationary point
 at most
$O(\epsilon^{-(p+1)/(p-q+1)})$ evaluations of the objective function and its first $p$
derivatives (whenever they exist).
Our model uses the underlying rotational symmetry of the Euclidean norm function to build a
  Lipschitzian approximation for the non-Lipschitzian group sparsity terms, which are defined by the group $\ell_2$-$\ell_a$ norm with $a\in (0,1)$.   The new
result shows that the partially-separable structure and non-Lipschitzian group sparsity terms in the objective
function may not affect the worst-case evaluation complexity order. }
\end{abstract}

{\small
\textbf{Keywords:} complexity theory,
nonlinear optimization,
non-Lipschitz functions,
partially-separable problems, group sparsity, isotropic model.
}

{\small
\textbf{AMS subject classifications:},  90C30, 90C46, 65K05

}

{\footnotesize \numsection{Introduction}}

Both applied mathematicians and computer scientists have, in recent years,
made significant contributions to the fast-growing field of worst-case
complexity analysis for nonconvex optimization (see \cite{CartGoulToin18a} for
a partial yet substantial bibliography). The purpose of this paper is to
extend the available general theory in two distinct directions.  The first is
to cover the case where the problem involving non-Lipschitzian group sparsity terms.
The second is to show that the ubiquitous
partially-separable structure (of which standard sparsity is a special case)
can be exploited without affecting the complexity bounds.

We consider the partially-separable convexly constrained nonlinear optimization
problem:
\begin{equation}
\label{problem}
 \min_{x \in \calF} f(x)
 = \sum_{i \in \calN} f_i(U_ix) + \sum_{i \in \calH} \|U_ix-b_i\|^a
 \eqdef \sum_{i \in \calN \cup \calH} f_i(U_ix)
\end{equation}
where
$\calN\cup\calH \eqdef \calM$,
$\calN\cap\calH = \emptyset$,
$f_i$ is a continuously $p$ times differentiable function from $\Re^{n_i}$
into $\Re$ for $i\in \calN$, and $f_i(x)=\|U_ix-b_i\|^a$ for $i\in \calH$,
$a \in (0,1)$,
$\|\cdot\|$ is the Euclidean norm,
$U_i\in \Re^{n_i \times n}$ with $n_i \leq n$,
and $b_i\in \Re^{n_i}$.
Without loss of generality, we assume that, for each $i \in \calM$, $U_i$ has
full row rank and $\|U_i\| = 1$, and that the ranges of the $U_i^T$ for $i\in
\calN$ span $\Re^n$ so that the intersection
of the nullspaces of the $U_i$ is reduced to the origin.
We also assume that the ranges of the $U_i^T$ (for $i
\in \calH$) are orthogonal, that is
\beqn{Ui-ortho-H}
U_iU_j^T = 0 \tim{for} i\neq j,\, i,j\in \calH.
\eeqn
Without loss of generality, we furthermore assume that the rows of $U_i$ are
orthonormal for $i \in \calH$. Our final assumption, as in
\cite{ChenToinWang17}, is that the feasible set  $\calF\subseteq \Re^n$ is
non-empty closed and convex, and that it is ``kernel-centered'' in that
is, if $P_\calX[\cdot]$ is the orthogonal projection onto the convex set
$\calX$ and $^{\dagger}$ denotes the Moore-Penrose generalized inverse,
then
\beqn{kernel-centered}
U_i^\dagger b_i + P_{\ker(U_i)}[\calF] \subseteq \calF
\tim{whenever} b_i \in U_i\calF, \quad i\in \calH.
\eeqn

These assumptions do not restrict our study for applications. For example, consider
 the row sparse problem in multivariate regression \cite{huang2009,huang2010benefit,obozinski2011}
 \begin{equation}\label{Prob4}
 \displaystyle  \min_{X\in R^{\nu\times \gamma}} \, \|HX-B\|^2_F +\lambda \|X\|_{\ell_a/\ell_2},
 \end{equation}
 where $H\in \Re^{\kappa\times \nu}, B \in \Re^{\kappa\times \gamma}$,$
 \|\cdot\|_F$
 is the Frobenius norm of a matrix,
 \[
 \|HX-B\|_F^2
 = \sum_{j=1}^{\gamma}\sum^\kappa_{i=1} (\sum_{\ell=1}^\nu H_{i\ell} X_{\ell, j}-B_{ij})^2
 \quad  {\rm and} \quad
 \|X\|_{\ell_a/\ell_2}=\sum^\nu_{i=1}\big(\sum^\gamma_{j=1} X_{ij}^2\big)^{\frac{a}{2}}.
 \]
 Let $n=\nu\gamma$, $\calF=\Re^n, b_i=0$,
 $x=(x_{11}, x_{12},\ldots, x_{\nu\gamma})^T\in \Re^n$ and set $U_i\in \Re^{\nu\times n}$ for $i\in \calN=\{1,\ldots, \gamma\}$
 be the projection whose entries are 0, or 1 such that $U_ix$ be the
 $i$th column of $X$ and $U_i\in \Re^{\gamma\times n}$ for $i\in \calH=\{1,\ldots, \nu\}$  be the
 projection whose entries are 0, or 1 such that $U_ix$ be the $i$th row of $X$. Then problem \req{Prob4} can be written in the form of
 (\ref{problem}).  It is easy to see that the $\{U_i^T\}_{i \in \calN}$  span $\Re^n$.   Hence, all assumptions mentioned above hold
 for problem \req{Prob4}.

Problem (\ref{problem}) encompasses the non-overlapping group sparse
optimization problems.  Let $G_1, \ldots, G_m$ be subsets of $\{1,\ldots, n\}$
representing known groupings of the decision variable with size $n_1,\ldots,
n_m$ and $G_i\cap G_j=\emptyset, i\neq j$.  In this case, problem
\req{problem} reduces to
\begin{equation}\label{Prob2}
\min_{x \in \calF} {\displaystyle f_1(x) + \lambda \sum^m_{i=1}\|U_ix\|^a},
 \end{equation}
where $f_1:\Re^n\to \Re_+$ is a smooth loss function, $\lambda >0$ is a positive number and
$U_i \in \Re^{n_i\times n}$ is defined in the following way
\[
(U_i)_{kj}=\left\{\begin{array}{ll}
1  &  \, {\rm if} \, j\in G_i\\
0  &  \, {\rm otherwise}
\end{array} \quad \quad {\rm for} \quad k=1,\ldots, n_i.
\right.
\]
Thus $U_ix=x_{G_i}$ is the  $i$th group variable vector in $\Re^{n_i}$ with components
$x_j, j\in G_i.$  If $\calF=\{x \, | \, \alpha_i \le x_i\le \beta_i, i=1,\ldots,n \}$
with $\alpha_i <0 <\beta_i$, then all assumptions mentioned above  with $U_1=I\in \Re^{n\times n}$,
for $1\in \calN$ hold for problem \req{Prob2}.

In problem (\ref{problem}), the decision variables have a group structure so
that components in the same group tend to vanish simultaneously. Group
sparse optimization problems have been extensively studied in recent years due
to numerous applications.  In machine learning and statistics, when the
explanatory variables have high correlative nature or can be naturally
grouped, it is important to study variable selection at the group sparsity
setting \cite{breheny2015,huang2009,huang2010benefit,Lee2016,MaHuang,Yuan}.
In compressed sensing, group sparsity is refereed to as block sparsity and has
been efficiently used to recovery signals with special block structures
\cite{Ahsen,Eldar,juditsky2012,Lee2012,Lv}.  In spherical harmonic
representations of random fields on the sphere, group Lasso penalty grouped
the coefficients of homogeneous harmonic polynomials of the same degree is
rotationally invariant while Lasso penalty ($G_i=\{i\}$) is not
\cite{gia2018}.

Problem (\ref{problem}) with $a\in (0,1)$ and $n_i=1, i\in \calH$ has been
studied in \cite{Bian_Chen_SIOPT,Bian-Chen-Ye,CXY,CNY,CGWY,Chen_Rob}.  Chen,
Toint and Wang \cite{ChenToinWang17} show that an adaptive regularization
algorithm using a $p$-th order Taylor model for $p$ \emph{odd} needs in
general at most $O(\epsilon^{-(p+1)/p})$ evaluations of the objective function
and its derivatives (at points where they are defined) to produce an
$\epsilon$-approximate first order critical point. Since this complexity bound
is identical in order to that already known for convexly constrained
Lipschitzian minimization, the result in \cite{ChenToinWang17} shows that
introducing non-Lipschitzian singularities in the objective function may not
affect the worst-case evaluation complexity order.

The unconstrained optimization of smooth partially-separable was first considered in Griewank and Toint \cite{GrieToin82a}, studied by many researchers
\cite{GoldWang93,Gay96,ChenDengZhan98,MareRichTaka14,ConnGoulSartToin96a,ConnGoulToin00}
and extensively used in the popular
{\sf CUTEst} testing environment \cite{
GoulOrbaToin15b} as well as in the AMPL \cite{FourGayKern87}, {\sf LANCELOT}
\cite{ConnGoulToin92} and {\sf FILTRANE} \cite{GoulToin07b}
packages.

In problem (\ref{problem}), all these ``element
functions'' $f_i$  depend on $U_ix \in \Re^{n_i}$ rather than on $x$, which is most useful when $n_i \ll n$.  Letting
$$x_i=U_ix \in \Re^{n_i} \,\,   {\rm for} \,\,  i\in \calM \quad  {\rm and} \quad  f_\calI(x)= \sum_{i\in
  \calI}f_i(x) \, \, \mbox{for any} \,\,  \calI \subseteq \calM,$$
   we define
\[
f_\calN(x) \eqdef \sum_{i \in \calN} f_i(U_ix) = \sum_{i \in \calN} f_i(x_i)
\quad \, {\rm and} \quad\,
f_\calH(x) \eqdef \sum_{i \in \calH} f_i(U_ix) = \sum_{i \in \calH} f_i(x_i).
\]
The $p$-th degree Taylor series
\beqn{taylor}
T_{f_\calN,p}(x,s) = f_\calN(x) + \sum_{j=1}^p \frac{1}{j!}\nabla_x^jf_\calN(x)[s]^j,
\tim{where}
\nabla_x^jf_\calN(x)[s]^j
= \sum_{i \in \calN} \nabla_{x_i}^jf_i(x_i)[U_is]^j,
\eeqn
indicates that, for each $j$, only the $|\calN|$ tensors
$\{\nabla_{x_i}^jf_i(x_i)\}_{i\in \calN}$ of dimension $n_i^j$ needs to be computed
and stored.  Exploiting derivative tensors of order larger than 2
--- and thus using the high-order Taylor series \req{taylor} as a local model
of $f_{\calN}(x+s)$ in the neighbourhood of $x$ --- may therefore be practically
feasible in our setting since $n_i^j$ is typically orders of magnitude
  smaller than $n$.  The same comment applies to
$f_\calH(x)$ whenever $\|U_ix-b_i\|\neq 0.$

The main contribution of this paper is twofold.
\begin{itemize}
\item We propose a partially separable adaptive regularization algorithm with
  a $p$-th order Taylor model which uses the underlying rotational symmetry
  of the Euclidean norm function for $f_\calH$ and the first $p$ derivatives
  (whenever they exist) of the ``element functions'' $f_i$, for $i \in \calM$.
\item We show that the algorithm can produce an
  $(\epsilon,\delta)$-approximate $q$-th-order critical point of problem
  \req{problem} at most $O(\epsilon^{-(p+1)/(p-q+1)})$ evaluations of the
  objective function and its first $p$ derivatives for any $q \in \ii{p}$.
\end{itemize}

Our results extend worst-case evaluation complexity bounds for smooth
nonconvex optimization in \cite{CartGoulToin18b,CartGoulToin18a} which
do not use the structure of partially separable functions and do not
consider the Lipschitzian singularity. Moreover, our results subsume
the results for non-Lipschitz nonconvex optimization in \cite{ChenToinWang17}
which only consider the complexity with $q=1$ and $n_i=1$ for $i\in \calH$.

This paper is organized as follows.  In Section 2, we define an $(\epsilon,
\delta)$ $q$-order necessary optimality conditions for local minimizers of
problem \req{problem}.  A Lipschitz continuous model to approximate $f$ is
proposed in Section~3. We then propose the partially separable adaptive
regularization algorithm using the $p$-th order Taylor model in Section 4. In
Section 5, we show that the algorithm produces an
$(\epsilon,\delta)$-approximate $q$-th-order critical point at most
$O(\epsilon^{-(p+1)/(p-q+1)})$ evaluations of $f$ and its first $p$
derivatives.


We end this section by introducing notations used in the next four sections.

\noindent
{\bf Notations.}
For a
symmetric tensor $S$ of order $p$,
$S[v]^p$ is the result of applying $S$ to $p$ copies of the vector $v$ and
\beqn{Tnorm}
\|S\|_{[p]} \eqdef \max_{\|v\|=1}  | S [v]^p |
\eeqn
is the associated induced norm for such tensors. If $S_1$ and $S_2$ are tensors,
$S_1\otimes S_2$ is their tensor product and $S_1^{k\otimes}$ is the product
of $S_1$ $k$ times with itself. For any set $\calX$, $|\calX|$ denotes its
cardinality.

Because the notion of partial separability hinges on geometric interpretation
of the problem, it is useful to introduce the various subspaces of interest
for our analysis.  We will extensively use the following definitions.
As will become clear in Section~\ref{optimality-s}, we will need to
identify
\beqn{Cxeps-def}
\calC(x,\epsilon) \eqdef \{ i \in \calH \mid \|U_ix-b_i\| \leq \epsilon \}
\tim{and}
\calA(x,\epsilon) \eqdef \calH \setminus \calC(x,\epsilon),
\eeqn
the collection of hard elements which are close to singularity for a given
$x$ and its complement (the ``active'' elements), and
\beqn{Rxeps-def}
\calR(x,\epsilon)
\eqdef \bigcap_{i\in\calC(x,\epsilon)} \ker(U_i)
= \left[\spanset_{i\in\calC(x,\epsilon)}(U_i^T)\right]^\perp
\eeqn
the subspace in which those nearly singular elements are invariant.
(When $\calC(x,\epsilon)=\emptyset$, we set $\calR(x,\epsilon)=\Re^n$.)
For convenience, if $\epsilon=0$, we denote $\calC(x)\eqdef \calC(x,0)$,
$\calA(x) \eqdef \calA(x,0)$, $\calR(x)\eqdef \calR(x,0)$ and
$\calW(x)\eqdef \calW(x,0)$.
From these definitions, we have
\beqn{dd}
U_id = 0,  \tim{for}  i\in \calC(x),  \,\,\,  d\in \calR(x).
\eeqn
Also denote by
\beqn{Rii-def}
\calR_{\{i\}} \eqdef \spanset(U_i^T)
\eeqn
and observe that \req{Ui-ortho-H} implies that the $\calR_{\{i\}}$ are
orthogonal for $i \in \calH$. Hence $\calR_{\{i\}}$ is also
the subspace in which all singular elements are invariant but the $i$-th.
We also denote the ``working'' collection of elements not close to singularity
by
\beqn{Wxeps-def}
\calW(x,\epsilon) \eqdef \calN \cup \calA(x,\epsilon).
\eeqn
 If $\{x_k\}$ is a sequence of iterates
in $\Re^n$, we also use the shorthands
\beqn{CRW-short}
\calC_k = \calC(x_k,\epsilon),
\ms
\calA_k = \calA(x_k,\epsilon),
\ms
\calR_k = \calR(x_k,\epsilon)
\tim{and}
\calW_k = \calW(x_k,\epsilon).
\eeqn
We will make frequent use of
\beqn{fWk-def}
f_{\calW_k}(x) \eqdef \sum_{i \in \calW_k}f_i(x),
\eeqn
which is objective function ``reduced'' to the elements ``away from singularity''
at $x_k$.

For some $x, s\in \Re^n$, we  often use the notations $
r_i = U_ix - b_i
$ and $s_i=Us$.

\numsection{Necessary optimality conditions}\label{optimality-s}

At variance with the theory developed in \cite{ChenToinWang17}, which solely
covers convergence to $\epsilon$-approximate first-order stationary points, we
now consider arbitrary orders of optimality. To this aim, we follow
\cite{CartGoulToin18b} and define, for a sufficiently smooth function
$h:\Re^n\rightarrow\Re$ and a convex set $\calF\subseteq\Re^n$, the vector $x$ to be an
$(\epsilon,\delta)$-approximate $q$-th-order stationary point ($\epsilon>0, \delta>0$, $q\in\ii{p}$)
of $\min_{x\in \calF} h(x)$ if, for some $\delta \in (0,1]$
\beqn{optimality-smooth}
\phi_{h,q}^\delta(x) \leq \epsilon \chi_q(\delta)
\eeqn
where
\beqn{phi-def}
\phi_{h,q}^\delta(x)
\eqdef h(x)-\min_{\stackrel{x+d\in \calF}{\|d\|\leq\delta}}T_{h,q}(x,d),
\eeqn
and
\beqn{chi-def}
\chi_q(\delta) \eqdef \sum_{\ell=1}^q \frac{\delta^\ell}{\ell!}.
\eeqn
In other words, we declare $x$ to be an $(\epsilon,\delta)$-approximate
$q$-th-order stationary point if the scaled maximal decrease that can be
obtained on the $q$-th order Taylor series for $h$ in a neighbourhood of $x$ of
radius $\delta$ is at most $\epsilon$.  We refer the reader to
\cite{CartGoulToin18b} for a detailed motivation and discussion of this
measure. For our present purpose, it is enough to observe that
$\phi_{h,q}^\delta(x)$ is a continuous function of $x$ and $\delta$ for any $q$.  Moreover,
for $q=1$ and $q=2$, $\delta$ can be chosen equal to one and $\phi_{h,1}^1(x)$
and $\phi_{h,2}^1(x)$ are easy to compute. In the unconstrained case,
\[
\phi_{h,1}^1(x) = \|\nabla_x^1h(x)\|
\]
and computing $\phi_{h,2}^1$ reduces to solving the standard trust-region problem
\[
\phi_{h,2}^1(x) = \left|\min_{\|d\|\leq 1} \nabla_x^1h(x)[d] + \half
\nabla_x^2h(x)[d]^2 \right|.
\]
In the constrained case,
\[
\phi_{h,1}^1(x) = \left| \min_{\mystack{x+d \in \calF}{\|d\|\leq
    1}}\nabla_x^1h(x)[d]\right|,
\]
which is the optimality measure used in \cite{CartGoulToin12b} or
\cite{ChenToinWang17} among others. However, given
the potential difficulty of solving the global optimization problem in
\req{phi-def} for $q>2$, our approach remains, for now, conceptual for such
high optimality orders.

We now claim that we can extend the definition \req{optimality-smooth} to
cover problem \req{problem} as well.  The key observation is that, by the
definition of $\calW(x,\epsilon)$ and $\calR(x,\epsilon)$,
\beqn{f=fWonReps}
f_{\calW(x,\epsilon)}(x)=f_{\calW(x,\epsilon)}(x+d)\le f(x+d)
\leq f_{\calW(x,\epsilon)}(x+d)+\epsilon^a|\calH| \tim{for all} d \in \calR(x,\epsilon).
\eeqn
Note now that $f_{\calW(x,\epsilon)}$ is smooth around $x$ because it only contains
elements which are away from singularity, and hence that
$T_{f_{\calW(x,\epsilon)},p}(x,s)$ is well-defined. We may therefore define
$x$ to be an $(\epsilon,\delta)$-approximate $q$-th-order stationary point
for \req{problem} if, for some $\delta \in (0,1]$
\beqn{optimality}
\psi_{f,q}^{\epsilon,\delta}(x) \leq \epsilon \chi_q(\delta),
\eeqn
where we define
\beqn{psi-def}
\psi_{f,q}^{\epsilon,\delta}(x)\eqdef
f(x)-\min_{\mystack{x+d\in\calF}
                      {\|d\|\leq\delta,\, d\in\calR(x,\epsilon)}} T_{f_{\calW(x,\epsilon)},q}(x,d).
\eeqn
By the definition of $\calW(x,\epsilon)$, we have $f_{\calW(x,\epsilon)}(x)\le f(x)$ and thus
\[f_{\calW(x,\epsilon)}(x)-\min_{\mystack{x+d\in\calF}
                      {\|d\|\leq\delta,\, d\in\calR(x,\epsilon)}} T_{f_{\calW(x,\epsilon)},q}(x,d)    \le  \psi_{f,q}^{\epsilon,\delta}(x).
\]

Taking $\epsilon=0$, $x$ is $q$-th-order stationary point if
$\psi_{f_{\calW(x)},q}^{0,\delta} =0$, as we now prove.

\begin{theorem}
If $x_*$ is a local minimizer of (\ref{problem}), then there is $\delta \in (0,1]$ such that
\begin{equation}\label{Op1}
\psi_{f,q}^{0,\delta}(x_*)=0.
\end{equation}
\end{theorem}
\proof{
Suppose first that $\calR(x_*) = \{0\}$ (which happens if there exists $x_* \in \calF$ such that $f_\calH(x_*)=0$ and
$\spanset_{i\in \calH}\{U_i^T\} = \Re^n$).  Then \req{Op1}
holds vacuously with any $\delta\in (0,1]$. Now suppose that $\calR(x_*) \neq \{0\}$.
Let
\[
\delta_1=\min\left[1, \min_{i\in \calA (x_*)} \|U_ix_*-b_i\|\right] \in (0,1].
\]
Since $x_*$ is a local minimizer of \req{problem}, there
exists $\delta_2>0$ such that
\begin{equation*}
\begin{split}
f(x_*) & = \bigmin_{\mystack{x_*+d \in \calF}{\|d\|\leq\delta_2}}f_{\calN}(x_*+d)+\sum_{i\in\calH}\|U_i(x_*+d)-b_i\|^a\\
       & \leq \bigmin_{\mystack{x_*+d \in \calF}{\|d\|\leq\delta_2,\,d \in \calR(x_*)}}f_{\calN}(x_*+d)+\sum_{i\in\calH}\|U_i(x_*+d)-b_i\|^a\\
       &   =  \bigmin_{\mystack{x_*+d \in \calF}{\|d\|\leq\delta_2,\,d \in \calR(x_*)}}f_{\calN}(x_*+d)+\sum_{i\in\calA(x_*)}\|U_i(x_*+d)-b_i\|^a\\
       &   =  \bigmin_{\mystack{x_*+d \in \calF}{\|d\|\leq\delta_2,\,d \in \calR(x_*)}}f_{\calW(x_*)}(x_*+d),
\end{split}
\end{equation*}
where we used (\ref{dd}) and (\ref{Wxeps-def}) to derive the last two equalities, respectively.

Now we consider the reduced problem
\beqn{reduced}
\min_{\mystack{x_*+d\in \calF}{\|d\|\leq\delta_2,\, d \in\calR(x_*)}}f_{\calW(x_*)}(x_*+d).
\eeqn
Since we have that
\[
f_{\calW(x_*)}(x_*)
=f_{\calN}(x_*)+\sum_{i\in\calA(x_*)}\|U_i x_*-b_i\|^a
=f_{\calN}(x_*)+\sum_{i\in\calH}\|U_i x_*-b_i\|^a
=f(x_*),
\]
we obtain that
\[
f_{\calW(x_*)}(x_*)
\leq \min_{\mystack{x_*+d \in \calF}{\|d\|\leq\delta_2,\,d \in \calR(x_*)}}f_{\calW(x_*)}(x_*+d)
\]
and $x_*$ is a local minimizer of problem
\req{reduced}.

Note that for any $x_*+d$ in the ball $B(x^*, \delta_3)$ with $\delta_3<\delta_1$, we have
\[
\|U_i(x_*+d)-b_i\|
\geq \|U_ix_*-b_i\|-\|U_id\|
\geq \delta_1-\|U_i\|\|d\|
= \delta_1-\delta_3 > 0, \ms i\in \calA(x_*).
\]
Hence $f_{\calW(x_*)}(x_*+d)$ is $q$-times continuously differentiable, and has
Lipschitz continuous derivatives of orders 1 to $q$ in $B(x^*, \delta_3)$. By
Theorem 3.1 in \cite{CartGoulToin17c}, there is a  $\delta \in \big(0,
\min[\delta_2, \delta_3]\big]$, such that
\[
\psi_{f_{\calW(x_*)},q}^{0,\delta}(x_*)
=f_{\calW(x_*)}(x_*)-\min_{\mystack{x_*+d\in \calF}{\|d\|\leq\delta,\ d\in \calR(x_*)}} T_{f_{\calW(x_*)},q}(x_*,d)=0
\]
This, together with
$f(x_*)=f_{\calW(x_*)}(x_*)$,
gives the desired result \req{Op1}.
}

We call $x_*$ is
a $q$-th-order stationary point of \req{problem} if there is $\delta \in (0,1]$ such that
(\ref{Op1}) holds.


\lthm{epsilon_q}{
  For each $k$, let $x_k$ be an $(\epsilon_k, \delta_k)$-approximate
  $q$-th-order stationary point of \req{problem} with
  $1 \geq \delta_k\geq \bar{\delta}>0$ and $\epsilon_k\rightarrow 0$.
  Then any cluster point of $\{x_k\}$ is
  a $q$-th-order stationary point of \req{problem}.
}

\proof{
Let $x_*$ be a cluster  point of $\{x_k\}$. Without loss of generality, we
assume that  $x_* = \lim_{k\rightarrow\infty}x_k$.
From  $0< \chi_q(\delta)\le 2$ and $\psi_{f,q}^{\epsilon, \delta}(x)\ge 0$, we
have from \req{optimality} that
$\lim_{k\rightarrow\infty}\psi_{f,q}^{\epsilon_k, \delta_k}(x_k)=0$. We now
need to prove that $\psi_{f,q}^{0, \bar{\delta}}(x_*)=0$.

If $\calR(x_*)=\{0\}$, \req{Op1} holds vacuously with any $\delta >0$, and
hence $x_*$ is a $q$th-order-necessary minimizer of \req{problem}.
Suppose now that $\calR(x_*) \neq \{0\}$.
We first claim that there exists a $k_*\geq 0$ such that
\beqn{CkCstar}
\calC(x_k,\epsilon_k ) \subseteq \calC(x_*)
\tim{ for } k \geq k_*.
\eeqn
To prove this inclusion, we choose $k_*$
sufficiently large to ensure that
\beqn{kstar-def}
\|x_k-x_*\| + \epsilon_k < \min_{j \in \calA(x_*)}\|U_jx_*-b_j\|,
\tim{ for } k \geq k_*.
\eeqn
Such  a $k_*$ must exist, since the right-hand side of this inequality is strictly positive
by definition of $\calA(x_*)$.
For an arbitrary $k \geq k_*$ and an index  $i \in
\calC(x_k,\epsilon_k )$, using the definition of $\calC(x,\epsilon)$, the
identity $\|U_i\| = 1$ and \req{kstar-def}, we  obtain that
\[
\|U_ix_*-b_i\|
\leq \|U_i(x_*-x_k)\| + \|U_ix_k -b_i\|
\leq \|x_* - x_k\| + \epsilon_k
< \min_{j \in  \calA(x_*)}\|U_jx_*-b_i\|.
\]
This implies that $\|U_ix_*-b_i\|= 0$ and $i \in \calC(x_*)$. Hence \req{CkCstar} holds.
By the definition of $\calR(x,\epsilon)$ and $\calW(x, \epsilon)$,
\req{CkCstar} implies that, for all $k$,
\beqn{RkRstar}
\calR(x_*) \subseteq \calR(x_k,\epsilon_k )
\tim{ and }
\calW(x_*) \subseteq \calW(x_k,\epsilon_k ).
\eeqn

For any fixed $k\ge k_*$,
consider  now the following three minimization problems:
\beqn{prob-a}
\hspace*{6mm}(A,k) \ms \left\{
\begin{array}{cl}
  \min_d      & \quad T_{f_{\calW(x_k,\epsilon_k)},q}(x_k,d)  \\
  \mbox{s.t.} & \quad x_k+d\in\calF,\,d\in\calR(x_k,\epsilon_k), \,\|d\|\le \delta_k,
\end{array}
\right.
\eeqn
\beqn{prob-b}
(B,k) \ms \left\{
\begin{array}{cl}
  \min_d      & \quad T_{f_{\calW(x_k,\epsilon_k)},q}(x_k, d)  \\
  \mbox{s.t.} & \quad x_k+d\in\calF,\, d\in\calR(x_*), \, \|d\|\le \delta_k,
\end{array}
\right.
\eeqn
and
\beqn{prob-c}
(C,k) \ms \left\{
\begin{array}{cl}
  \min_d      & \quad T_{f_{\calW(x_*)},q}(x_k, d)  \\
  \mbox{s.t.} & \quad x_k+d\in\calF, \, d\in\calR(x_*), \, \|d\|\le \delta_k.
\end{array}
\right.
\eeqn
Since $d=0$ is a feasible point of these three problems, their minimum values,
which we respectively denote by $\vartheta_{A,k}$, $\vartheta_{B,k}$ and
$\vartheta_{C,k}$, are all smaller than $f(x_k)$. Moreover, it follows
from the first part of \req{RkRstar} that, for each $k$,
\beqn{vartBA}
\vartheta_{B,k} \geq \vartheta_{A,k}.
\eeqn
It also follows from \req{CkCstar} and \req{Ui-ortho-H} that
\[
T_{f_{\calW(x_k,\epsilon_k)},q}(x_k, d)=
T_{f_{\calW(x_*)},q}(x_k, d) - f_{\calW(x_*)}(x_k) + f_{\calW(x_k,\epsilon_k)}(x_k)
\leq T_{f_{\calW(x_*)},q}(x_k, d) + |\calH|\epsilon_k^a
\]
for all $d\in \calR(x_*)$, and thus \req{vartBA} becomes
\beqn{vartABC}
\vartheta_{A,k} \leq \vartheta_{B,k} \le \vartheta_{C,k}+|\calH|\epsilon_k^a
\tim{ for all}k\ge k_*.
\eeqn
The assumption that  $x_k$ is an $(\epsilon_k, \delta_k)$-approximate
$q$th-order necessary minimizer of \req{problem} implies that
\beqn{vart2}
0 \leq f(x_k)-\vartheta_{C,k}-|\calH|\epsilon_k^a
\leq f(x_k)-\vartheta_{A,k}\le\epsilon_k \chi_q(\delta_k),
\tim{for all}  k \geq k_*.
\eeqn
Now (\ref{RkRstar}) implies that
$T_{f_{ \calW(x_*)},q}(x_k,d)\leq T_{f_{\calW(x_k)},q}(x_k,d)\leq T_{f,q}(x_k,d)$.  Hence
\[
f(x_k)-\vartheta_{C,k}=f(x_k)-\min_{\mystack{x_k+d\in
    \calF}{\|d\|\le \delta_k, \ d\in \calR(x_*)}}T_{f_{ \calW(x_*)},q}(x_k,d)
\ge\phi_{f,q}^{0,\delta_k}(x_k).
\]
As a consequence, \req{vart2} implies that
\beqn{epsilonc}
\phi_{f,q}^{0,\delta_k}(x_k)
\leq \epsilon_k \chi_q(\delta_k)+|\calH|\epsilon_k^a.
\eeqn
In addition, the feasible sets of the three problems \req{prob-a}-\req{prob-c}
are convex, and the objectives functions are polynomials with degree $q$.
By the perturbation theory for optimization problems
\cite[Theorem~3.2.8]{ConnGoulToin00}, we can claim that
\beqn{chiCk}
\lim_{k \rightarrow \infty} \vartheta_{C,k}
=\min_{\mystack{x_*+d\in \calF}{\|d\|\le \delta_*, \ d\in \calR(x_*)}}T_{f_{ \calW(x_*)},q}(x_*,d),
\eeqn
where $\delta_*=\liminf_{k \rightarrow \infty}\delta_k\geq \bar{\delta}$.
This implies that letting $k\to \infty$ in \req{epsilonc} gives
\[
\psi_{f,q}^{0, \bar{\delta}}(x_*)
= f(x_*)
  -\min_{\mystack{x_*+d\in \calF}{\|d\|\leq \bar{\delta},\,d\in\calR(x_*)}} T_{f_{ \calW(x_*)},q}(x_*,d)=0.
\]
}

We conclude this section by an important observation.  The optimality measure
\req{optimality} may give the impression (in particular in its use of
$\calR(x,\epsilon)$) that the ``singular'' and ``smooth''
parts of the problem are merely separated, and that one could possibly apply the
existing theory for smooth problems to the latter.  Unfortunately, this
is not true, because the ``separation'' implied by \req{optimality} does
depend on $\epsilon$, and one therefore need to show that the complexity of
minimizing the ``non-singular'' part does not explode (in particular with the
unbounded growth of the Lispchitz constant) when $\epsilon$ tends to
zero. Designing an suitable algorithm and proving an associated complexity
result comparable to what is known for smooth problems is the main challenge
in what follows.

\numsection{A Lipschitz continuous model of $f_{\calW_k}(x+s)$}\label{model-s}

Our minimization algorithm, described in the next section, involves the
approximate minimization of a \emph{model} $m(x_k,s)$ of $f_{\calW_k}$ in the
intersection of a neighbourhood of $x_k$ with $\calR_k$.  This model,
depending on function and derivatives values computed at $x_k$, should be
able to predict values and derivatives of $f$ at some neighbouring point
$x_k+s$ reasonably accurately. This is potentially difficult if the current
point happens to be near a singularity.

Before describing our proposal, we need to state a useful technical result.

\llem{dersnorm-l}{Let $a$ be a positive number and $r\neq 0$.
  Define, for a positive integer $j$,
  \beqn{factor-def}
  \pi(a-j) \eqdef a\prod_{i=1}^{j-1}(a-i).
  \eeqn
  Then, if $\nabla_\cdot^j\big\|r\big\|^a$ is the value
  of the $j$-th derivative tensor of the function $\|\cdot\|^a$ with
  respect to its argument, evaluated at $r$, we have that,
  \beqn{ders}
  \nabla_\cdot^j \|r\|^a
   = \bigsum_{i=1}^j \phi_{i,j}
   \|r\|^{a-2i} \, r^{(2i-j) \otimes} \otimes I^{(j-i)\otimes}
  \eeqn
  for some scalars $\{\phi_{i,j}\}_{i=1}^j$  such that
  $\sum_{i=1}^j \phi_{i,j} = \pi(a-j)$ , and that
  \beqn{dersnorm}
  \big\|\,\nabla_\cdot^j \|r\|^a \,\big\|_{[j]} = |\pi(a-j)|\,\|r\|^{a-j}.
  \eeqn
  Moreover, if $\beta_1, \beta_2$ are positive reals and $\|r\|=1$, then
  \beqn{difdersnorm}
  \left\|\,\nabla_\cdot^j \|\beta_1 r\|^a - \nabla_\cdot^j \|\beta_2 r\|^a\,\right\|_{[j]}
  = |\pi(a-j)|\,\left|\beta_1^{a-j} - \beta_2^{a-j}\right|.
  \eeqn
}

\proof{See appendix. }

Consider now the elements $f_i$ for $i\in \calN$. Each such element is $p$ times
continuously differentiable and, if we assume that its $p$-th derivative
tensor $\nabla_x^pf_i$ is globally Lipschitz continuous with constant $L_i
\geq 0$ in the sense that, for all $x_i,y_i \in \Re^{n_i}$
\beqn{tensor-Lip-fi}
\|\nabla_{x_i}^pf_i(x_i) - \nabla_{x_i}^p f_i(y_i)\|_{[p]} \leq L_i \|x_i-y_i\|,
\eeqn
then it can be shown (see \cite[Lemma~2.1]{CartGoulToin18b}) that
\beqn{f-Lip-1}
f_i(x_i+s_i) = T_{f_i,p}(x_i,s_i) + \frac{1}{(p+1)!} \tau_i L_i \|s_i\|^{p+1}
\tim{ with }| \tau_i | \leq 1.
\eeqn
Because $\tau_i L_i$ in \req{f-Lip-1} is usually unknown in
practice, it may not be possible to use \req{f-Lip-1} directly to model $f_i$
in a neighbourhood of $x$.  However, we may replace this term with an
adaptive parameter $\sigma_i$, which yields the following
$(p+1)$-rst order model for the $i$-th ``nice'' element
\beqn{miN-def}
m_i(x_i,s_i)=T_{f_i,p}(x_i,s_i)+\frac{1}{(p+1)!}\ \sigma_i\|s_i\|^{p+1},
\ms  (i \in \calN).
\eeqn
Using the associated Taylor's expansion
would indeed ignore the non-Lipschitzian singularity occurring for $r_i=0$
and this would restrict the validity of the model to a possibly very small
neighbourhood of $x_k$ whenever $r_i$ is small for some $i \in
\calA(x_k,\epsilon)$.  Our proposal is to use the underlying rotational
symmetry of the Euclidean norm function to build a better Lipschtzian
model. Suppose that $r_i \neq 0 \neq r_i+s_i$ and let
\beqn{uiui+-def}
u_i = \frac{r_i}{\|r_i\|},
\ms
r_i^+ = r_i(x+s) = r_i+s_i
\tim{ and }
u_i^+ = \frac{r_i^+}{\|r_i^+\|}.
\eeqn
Moreover, let $R_i$ be the rotation in the $(u_i,u_i^+)$ plane\footnote{If
  $u_i=u_i^+$, $R_i= I$. If $n_i=1$ and $r_ir_i^+<0$, this rotation is just
  the mapping from $\Re_+$ to $\Re_-$, defined by a simple sign change, as in
  the two-sided model of \cite{ChenToinWang17}.} such that
\beqn{roti-def}
R_iu_i^+ = u_i.
\eeqn
We observe that, given the isotropic nature of the Euclidean norm, the value
$\|r_i\|^a$ and of its derivatives with respect to $s_i$ can be deduced
from those $\| \, \|r_i\|u_i^+ \,\|^a$. More precisely, for any $d\in\Re^{n_i}$,
\beqn{rotation}
\big\|\, \|r_i\|u_i^+ \, \big\|^a = \|r_i\|^a
\tim{and}
\nabla_\cdot^\ell\big\| \,\|r_i\|u_i^+\,\big\|^a[d]^\ell
=\nabla_\cdot^\ell\|r_i\|^a[R_id]^\ell.
\eeqn
For example, when $j=1$, one
verifies that
\[
\begin{array}{ll}
\nabla_\cdot^1\|r_i\|^a[R_id]
& = a \|r_i\|^{a-2}r_i^TR_id \\*[1.5ex]
& = a \|r_i\|^{a-2} \|r_i\|( R_i^Tu_i)^Td \\*[1.5ex]
& = a \big\|\,\|r_i\| u_i^+\,\big\|^{a-2} (\|r_i\|u_i^+)^Td\\*[1.5ex]
& = \nabla_\cdot^1\big\| \,\|r_i\|u_i^+\,\big\|^a [d].
\end{array}
\]
We may then choose to compute the Taylor's expansion for the function $\|\cdot\|^a$
around $\|r_i\|u_i^+$, that is
\[
\begin{array}{ll}
\|r_i^+\|^a
& = \big\|\,\|r_i^+\|u_i^+\,\big\|^a \\*[1.5ex]
& =
\big\|\,\|r_i\|u_i^+\,\big\|^a+\bigsum_{\ell=1}^\infty\bigfrac{1}{\ell!}
     \nabla_\cdot^\ell\big\|\|r_i\|u_i^+\big\|^a
      \big[(\|r_i^+\|-\|r_i\|)u_i^+\big]^\ell\\*[1.5ex]
& = \|r_i\|^a
      +\bigsum_{\ell=1}^\infty\bigfrac{(\|r_i^+\|-\|r_i\|)^\ell}{\ell!}
       \nabla_\cdot^\ell\big\|r_i\big\|^a\big[R_iu_i^+\big]^\ell\\*[1.5ex]
& = \|r_i\|^a
      +\bigsum_{\ell=1}^\infty
       \bigfrac{(\|r_i^+\|-\|r_i\|)^\ell}{\ell!}
       \nabla_\cdot^\ell\big\|r_i\big\|^a[u_i]^\ell.
\end{array}
\]
Using the expression \req{ders} applied to $\ell$ copies of the unit vector
$u_i$ and the fact that $r_i^{j\otimes}[u_i]^j= (r_i^Tu_i)^j=\|r_i\|^j$ for all
$j\in\Na$, we now deduce that, for $\zeta_i= \|r_i+s_i\|-\|r_i\|\geq -\|r_i\|$,
\beqn{Taylorsi}
\|r_i+s_i\|^a
= \|r_i\|^a  + \bigsum_{\ell=1}^\infty
       \bigfrac{\pi(a-\ell)}{\ell!}\zeta_i^\ell
       \|r_i\|^{a-\ell},
\eeqn
which is nothing else than the Taylor's expansion of $\|r_i+\zeta_iu_i\|^a$
(or, equivalently, of $\big\| \|r_i\|u_i^+ + \zeta_iu_i^+\big\|^a$)
expressed as a function of the scalar variable $\zeta_i \geq -\|r_i\|$. As can be
expected from the isotropic nature of the Euclidean norm, the value of
$\|r_i^+\|^a$ (and of its derivatives after a suitable rotation) only
depend(s) on the distance of $r_i^+$ to the singularity at zero.
Thus, limiting the development \req{Taylorsi} to degree $p$ (as in
\req{miN-def}), it is natural to define
\beqn{TiH-def}
\mu(\|r_i\|, \zeta_i)
\eqdef \|r_i\|^a
      +\bigsum_{\ell=1}^p\bigfrac{\pi(a-\ell)}{\ell!}\zeta_i^\ell
      \|r_i\|^{a-\ell},\\*[1.5ex]
\ms
(i \in \calA(x,\epsilon)),
\eeqn
which is a unidimensional model of $\|r_i+\zeta_iu_i\|^a$ based of the
residual value $\|r_i\|$.
Note that $\mu(\|r_i\|,\zeta_i)$ is Lipschitz continuous as a function of
$\zeta_i$ as long as $\|r_i\|$ remains uniformly  bounded
away from zero, with a Lipschitz constant depending on the lower bound.  We
then define the isotropic model
\beqn{miH-def}
m_i(x_i,s_i) \eqdef \mu(\|r_i\|,\zeta_i) = \mu_i(\|r_i\|,\|r_i+s_i\|-\|r_i\|)
\tim{ for } i\in \calA(x,\epsilon),
\eeqn
so that, abusing notation slightly,
\[
m_i(x_i,s_i) = T_{m_i,p}(x_i,s_i) \eqdef \|r_i\|^a
      +\bigsum_{\ell=1}^p\bigfrac{\pi(a-\ell)}{\ell!}\zeta_i^\ell
      \|r_i\|^{a-\ell},\\*[1.5ex]
\ms
(i \in \calA(x,\epsilon)).
\]
\noindent
We now state some useful properties of this model.

\llem{miprop-l}{Suppose that $p$ is odd and that $\calA(x) \neq \emptyset$ for
  some $x \in \calF$.
  Then, for $i \in \calA(x)$,
  \beqn{overestH}
  m_i(x_i,s_i) \geq \|r_i+s_i\|^a,
  \eeqn
  and, whenever $\|r_i+s_i\|\leq \|r_i(x)\|$,
  \beqn{dersmi-odd}
  \nabla_\zeta^\ell \mu(\|r_i\|,\zeta_i))
  \geq \nabla_\zeta^\ell \mu(\|r_i\|,0)
  = \pi(a-\ell)\|r_i\|^{a-\ell}
  \ms
  (\ell \mbox{~odd}),
  \eeqn
  and
  \beqn{dersmi-even}
  \nabla_\zeta^\ell \mu(\|r_i|,\zeta_i)
  \leq \nabla_\zeta^\ell\mu(\|r_i\|,0)
  = \pi(a-\ell)\|r_i\|^{a-\ell}
  \ms
  (\ell \mbox{~even}).
  \eeqn
  As a consequence, $m_i(x_i,\zeta_i)$ is a concave function of $\zeta_i$ on
  the interval $[-\|r_i\|,0]$.
}

\proof{Let $i\in \calA(x)$.
From the mean-value theorem and \req{Taylorsi}, we have that, for some $\nu
\in (0,1)$,
\beqn{mvmi}
\|r_i+s_i\|^a
= \|r_i\|^a + \bigsum_{\ell=1}^p
       \bigfrac{\pi(a-\ell)}{\ell!}\zeta_i^\ell
       \|r_i\|^{a-\ell} + \bigfrac{\pi(a-p-1)}{(p+1)!} \zeta_i^{p+1}\|r_i+\nu
       \zeta_i u_i\|^{a-p-1}.
\eeqn
Since $p$ is odd, we obtain that $\pi(a-p-1)<0$ and $\zeta_i^{p+1} \geq 0$.
Thus \req{overestH} directly results from \req{mvmi},  \req{TiH-def}       and \req{miH-def}.
Now \req{TiH-def} and \req{miH-def} together imply that
\beqn{derzeta}
\nabla_\zeta^\ell \mu(\|r_i\|,\zeta_i)
= \nabla_\zeta^\ell \mu(\|r_i\|,0)
   + \sum_{j=\ell+1}^p\frac{\pi(a-j)}{(j-\ell)!}\zeta_i^{j-\ell}\|r_i\|^{a-j}
\eeqn
for $\zeta_i = \|r_i(x)+s_i\|-\|r_i(x)\| \leq 0$.
Suppose first that $\ell$ is odd.  Then we have that $\pi(a-j)$ is negative for
even $j$, that is exactly when $\zeta_i^{j-\ell}$ is non-positive.  Hence
every term in the sum of the right-hand side of \req{derzeta} is non-negative
and \req{dersmi-odd} follows.  Suppose now that $\ell$ is even. Then $\pi(a-j)$
is negative for odd $j$, which is exactly when $\zeta_i^{j-\ell}$ is
non-negative. Hence every term in the sum of the right-hand side of
\req{derzeta} is non-positive and \req{dersmi-even} follows. The last
conclusion of the lemma is then deduced by considering $\ell=2$ in
\req{dersmi-even} and observing that
$\pi(a-\ell)\|r_i(x)\|^{a-\ell}=a(a-1)\|r_i(x)\|^{a-2}<0$.
} 

\noindent
Thus the isotropic model $m_i(x_i,s_i)$ overestimates the true function
$\|r_i(x)+s_i\|^a$ and correctly reflects its concavity in the direction of
its singularity.  But $m_i(x_i,s_i)= \mu(\|r_i\|,\zeta_i)$ is now Lipschitz
continuous as a function of $s_i$, while $\|r_i(x)+s_i\|^a$ is not.

Combining \req{miN-def} and \req{miH-def} now allows us to define a model for
the complete $f$ on $\calR(x,\epsilon)$ as
\beqn{mfull-def}
m(x,s) 
\eqdef \sum_{i\in \calW(x,\epsilon)} m_i(x_i,s_i).
\eeqn
Since $r_i(x)\neq 0$ for $i \in \calA(x,\epsilon)$, this model is in turn
well defined.

We conclude this section by observing that writing the problem in the
partially-separable form \req{problem} is the key to expose the singular parts
of the objective function, which then allows exploiting their rotational
symmetry.

\numsection{The adaptive regularization algorithm}\label{algo-s}

Having defined a model of $f$, we may then use this model within a
regularization minimization method inspired from the AR$p$ algorithm in
\cite{CartGoulToin18b}. In such an algorithm, the step from an iterate $x_k$
is obtained by attempting to (approximately) minimize the model (\req{mfull-def}
in our case).  If an $(\epsilon,\delta)$-approximate $q$-th-order-necessary
minimizer is sought, this minimization is terminated as soon as the step $s_k$
is long enough, in that
\beqn{longs}
\|s_k\|\geq \varpi \epsilon^{\frac{1}{p-q+1}}
\eeqn
for some constant $\varpi \in (0,1]$,
or as soon as the trial point $x_k+s_k$ is an
approximate $q$-th-order-necessary minimizer \emph{of the model}, in the sense
that
\beqn{term-m}
\psi_{m,q}^{\epsilon,\delta_k}(x_k,s_k)
\leq \min\left[\frac{\theta\|s_k\|^{p-q+1}}{(p-q+1)!},
               \, a \min_{i\in \calA(x_k+s_k, \epsilon)}\|r_i(x_k+s_k)\| \right]\chi_q(\delta_k)
\eeqn
for some $\theta, \delta_k \in (0,1]$, where $\psi_{m,q}^{\epsilon,\delta_k}(x_k+s_k)$
is the optimality measure \req{psi-def} computed for the model $m(x,s)$,
that is
\beqn{psim-def}
\psi_{m,q}^{\epsilon,\delta}(x,s)
\eqdef m(x,s)-\min_{\mystack{x+s+d \in \calF}{\|d\|\leq\delta, d\in \calR(x,\epsilon)}}T_{m,q}(x,s+d).
\eeqn
In view of the optimality condition \req{optimality}, we also require that, if
$\|r_i(x+s)\| \leq \epsilon$ occurs for some $i \in \calH$ in the course of
the model minimization, the value of $r_i(x+s)$ is fixed, implying that the
remaining minimization is carried out on $\calR(x_k+s,\epsilon)$.  As a
consequence, the dimension of $\calR(x_k+s,\epsilon)$ (and thus of $\calR_k$) is
monotonically non-increasing during the step computation and across all
iterations.  It was shown in \cite[Lemma~2.5]{CartGoulToin18b} that, unless
$x_k$ is an $(\epsilon,1)$-approximate $p$-th-order-necessary minimizer (which
is obviously enough for the whole algorithm to terminate), a step satisfying
\req{term-m} can always be found.  The fact that this condition must hold on a
subspace of potentially diminishing dimension clearly does not affect the
result, and indicates that \req{term-m} is well-defined. This model
minimization is in principle simpler than the original problem because the
general nonlinear $f_i$ have been replaced by locally accurate polynomial
approximations and also because the model is now Lipschitz continuous, albeit
still non-smooth. Importantly, the model minimization \emph{does not involve
any evaluation of the objective function} or its derivatives, and model
evaluations within this calculation therefore do not affect the overall
evaluation complexity of the  algorithm.

We now introduce some useful notation for describing our algorithm. Define
\[
x_{i,k} \eqdef U_ix_k,
\ms r_{i,k} \eqdef U_ix_k-b_i,
\ms s_{i,k} \eqdef U_is_k,
\ms u_{i,k} \eqdef \frac{r_{i,k}}{\|r_{i,k}\|}
\]
and
\[
\calA_k^+ \eqdef \calA(x_k+s_k,\epsilon),
\ms
\calR_k^+ \eqdef \calR(x_k+s_k,\epsilon)
\tim{and}
\calW_k^+ \eqdef \calW(x_k+s_k,\epsilon).
\]
Also let
\[
\Delta f_{i,k} \eqdef f_i(x_{i,k}) - f_i(x_{i,k}+s_{i,k}),
\ms
\Delta f_k
\eqdef f_{\calW_k^+}(x_k) - f_{\calW_k^+}(x_k+s_k)
= \sum_{i \in \calW_k^+} \Delta f_{i,k},
\]
\[
\Delta m_{i,k} \eqdef m_i(x_{i,k},0) - m_i(x_{i,k},s_{i,k}),
\ms
\Delta m_k
= \sum_{i \in \calW_k^+} \Delta m_{i,k},
\]
and
\beqn{deltaT-def}
\begin{array}{lcl}
\Delta T_k
& \eqdef &  T_{f_{\calW_k^+},p}(x_k,0) - T_{f_{\calW_k^+},p}(x_k,s_k) \\*[2ex]
&    =   & [T_{f_\calN,p}(x_k,0) - T_{f_\calN,p}(x_k,s_k)]
+ \big[m_{\calA_k^+}(x_k,0) - m_{\calA_k^+}(x_k,s_k)\big]\\*[2ex]
&    =   & \Delta m_k+\bigfrac{1}{(p+1)!}\bigsum_{i \in \calN}\sigma_{i,k}\|s_{i,k}\|^{p+1}.
\end{array}
\eeqn
Our partially-separable adaptive regularization degree-$p$ algorithm PSAR$p$ is then
given by Algorithm~\ref{psarp} \vpageref{psarp}.

\algo{psarp}{Partially-Separable Adaptive Regularization (PSAR$p$)}{
\vspace*{-3mm}
\begin{description}
\item[Step 0: Initialization:] $x_0\in\calF$ and $\{\sigma_{i,0}\}_{i\in\calN}
  >0$ are given as well as the accuracy $\epsilon \in (0,1]$ and constants $0
<\gamma_0 < 1 <   \gamma_1 \leq \gamma_2$, $\eta \in (0,1)$, $\theta \geq0$, $\delta_{-1}=1$,
$\sigma_{\min} \in  (0, \min_{i\in\calN}\sigma_{i,0}]$ and  $\kap{big}>1$. Set $k =0$.
\item[Step 1: Termination:] Evaluate $f(x_k)$ and
  $\{\nabla_x^j f_{\calW_k}(x_k)\}_{j=1}^q$. If
\beqn{optimk}
\psi_{f,q}^{\epsilon,\delta_{k-1}}(x_k) \leq \epsilon \chi_q(\delta_{k-1})
\eeqn
  return $x_\epsilon=x_k$ and terminate. Otherwise evaluate
  $\{\nabla_x^j f_{\calW_k}(x_k)\}_{j=q+1}^p$.
\item[Step 2: Step computation:] Attempt to compute a
  step $s_k\in \calR_k$ such that $x_k+s_k \in \calF$, $m(x_k,s_k)<  m(x_k,0)$
  and either \req{longs} holds or \req{term-m} holds for some $\delta_k\in (0,1]$.
  If no such step exists, return $x_\epsilon=x_k$ and terminate.
\item[Step 3: Step acceptance:] Compute
\beqn{rhok-def}
\rho_k = \bigfrac{\Delta f_k}{\Delta T_k}
\eeqn
and set $x_{k+1} = x_k$ if $\rho_k < \eta$, or $x_{k+1} = x_k+s_k$ if $\rho_k \geq \eta$.

\item[Step 4: Update the ``nice'' regularization parameters:] For
$i \in \calN$, if
\beqn{up-cond}
f_i(x_{i,k}+s_{i,k}) > m_i(x_{i,k},s_{i,k})
\eeqn
set
\beqn{sig-incr}
\sigma_{i,k+1} \in [ \gamma_1 \sigma_{i,k}, \gamma_2 \sigma_{i,k} ].
\eeqn
Otherwise, if either
\beqn{muchover-neg}
\rho_k \geq \eta
\tim{ and }
\Delta f_{i,k} \leq 0
\tim{ and }
\Delta f_{i,k} < \Delta m_{i,k} - \kap{big}|\Delta f_k|
\eeqn
or
\beqn{muchover-pos}
\rho_k \geq \eta
\tim{ and }
\Delta f_{i,k} > 0
\tim{ and }
\Delta f_{i,k} > \Delta m_{i,k} + \kap{big}|\Delta f_k|
\eeqn
then set
\beqn{sig-decr}
\sigma_{i,k+1} \in [\max[\sigma_{\min}, \gamma_0\sigma_{i,k}], \sigma_{i,k} ],
\eeqn
else set
\beqn{sig-keep}
\sigma_{i,k+1} = \sigma_{i,k}.
\eeqn
Increment $k$ by one and go to Step~1.
\end{description}
}

Note that an $x_0 \in \calF$ can always be computed by projecting an
infeasible starting point onto $\calF$.  The motivation for the second and
third parts of \req{muchover-neg} and \req{muchover-pos} is to identify cases
where the isotropic model $m_i$ overestimates the element function $f_i$ to an
excessive extent, leaving some room for reducing the regularization and hence
allowing longer steps.  The requirement that $\rho_k \geq \eta$ in both
\req{muchover-neg} and \req{muchover-pos} is intended to prevent a situation
where a particular regularization parameter is increased and another decreased
at a given unsuccessful iteration, followed by the opposite situation at the
next iteration, potentially leading to cycling.

It is worthwhile to note the differences between the PSAR$p$ algorithm and
the algorithm discussed in \cite{ChenToinWang17}.
The first and most important is that the new algorithm is intended to find an
$(\epsilon,\delta)$-approximate $q$-th-order necessary minimizer for problem
\req{problem}, rather than a first-order critical point. This is made possible
by using the $q$-th-order termination criterion \req{optimk} instead of a
criterion only involving the first-order model decrease,  and by simultaneously
using the step termination criteria \req{longs} and \req{term-m} which again
replace a simpler version again based solely on first-order information.
The second is that the PSAR$p$ algorithm applies to the more general problem
\req{problem}, in particular using the isotropic model \req{TiH-def} to allow
$n_i>1$ for $i\in \calH$.

As alluded to above and discussed in \cite{CartGoulToin18a} and
\cite{BellGuriMoriToin18}, the potential termination of the algorithm in
Step~2 can only happen whenever $q>2$ and $x_k=x_\epsilon$ is an
$(\epsilon,1)$-approximate $p$-th-order-necessary minimizer within $\calR_k$,
which, together with \req{optimality}, implies that the same property holds
for problem \req{problem}. This is a significantly stronger optimality
condition than what is required by \req{optimk}. Also note that the
potentially costly calculation of \req{term-m} may be avoided if \req{longs}
holds.

Let the index set of the ``successful'' and ``unsuccessful'' iterations be
given by
\[
\calS \eqdef \{ k \geq 0 \mid \rho_k \geq \eta \}
\tim{and}
\calU \eqdef \{ k \geq 0 \mid \rho_k < \eta \}.
\]
Also define
\[
\calS_k \eqdef \calS \cap \iibe{0}{k}
\tim{ and }
\calU_k \eqdef \iibe{0}{k} \setminus \calS_k.
\]
We then state a bound on $|\calU_k|$ as a function of $|\calS_k|$. This is a
standard result for non-partially-separable problems (see
\cite[Theorem~2.4]{BirgGardMartSantToin17} for instance), but needs careful
handling of the model's overestimation properties to apply to our present
context.

\llem{a-succ-unsucc}{
Suppose that AS.2 and AS.3 hold and that $\sigma_{i,k}\leq \sigma_{\max}$ for
all $i \in \calM$ and all $k\geq0$.  Then, for all $k \geq 0$,
\[
k \leq \kappa^a |\calS_k| + \kappa^b,
\]
where
\[
\kappa^a \eqdef  1 + \frac{|\calN|\,|\log \gamma_0|}{\log \gamma_1}
\tim{ and }
\kappa^b \eqdef \frac{|\calN|}{\log \gamma_1}
                \log\left(\frac{\sigma_{\max}}{\sigma_{\min}}\right).
\]
}

\proof{See \cite[Lemma~4.11]{ChenToinWang17}. The proof hinges on \req{overestH}.}

\numsection{Evaluation complexity analysis}\label{complexity-s}

We are now ready for a formal analysis of the evaluation complexity of the
PSAR$p$ algorithm for problem \req{problem}, under the following assumptions.

\ass{AS.1}{ The feasible set $\calF$ is closed, convex, non-empty and
  kernel-centered (in the sense of \req{kernel-centered}).}
\vspace*{-4mm}
\ass{AS.2}{ Each element function $f_i$ ($i \in \calN$) is $p$ times
  continuously differentiable in an open set containing $\calF$, where $p$ is
  odd whenever $\calH \neq \emptyset$.}
\vspace*{-4mm}
\ass{AS.3}{The $p$-th derivative of each $f_i$ ($i \in \calN$) is
  Lipschitz continuous on $\calF$ with associated Lipschitz constant $L_i$ (in
  the sense of \req{tensor-Lip-fi}).
}
\vspace*{-4mm}
\ass{AS.4}{There exists a constant $f_{\rm low}$ such that $f_\calN(x) \geq f_{\rm
    low}$ for all $x \in \calF$.
}
\vspace*{-4mm}
\ass{AS.5}{If $\calH\neq\emptyset$, there exists a constant
  $\kappa_\calN \geq 0$ such that
  $ \|\nabla_x^j f_i(U_ix)\| \leq \kappa_\calN $ for all
  $x \in \calV$, $i \in \calN$ and $j \in \ii{p}$, where
  \beqn{calV-def}
  \calV \eqdef \left\{ x \in \calF \mid \tim{there exists} i \in \calH
  \tim{with}\|r_i(x)\| \leq \frac{a}{16}\right\}.
  \eeqn
}

\noindent
Note that AS.4 is necessary for problem \req{problem} to be well-defined.
Also note that, because of AS.2, AS.5 automatically holds if $\calF$ is
bounded or if the iterates $\{x_k\}$ remain in a bounded set. It is possible
to weaken AS.2 and AS.3 by replacing $\calF$ with the level set $\calL = \{ x
\in \calF \mid f(x) \leq f(x_0) \}$ without affecting the results
below. Finally observe that $\calV$ need not to be bounded, in particular if
$\spanset_{i\in \calH}(U_i)$ is a proper subspace of $\Re^n$. AS.5 is, of
course, unnecessary if $\calF$ or $\calV$ are bounded or the iterates remain
in a bounded set. The motivation for the particular choice of $\sfrac{1}{16}a$
in \req{calV-def} will become clear in Lemma~\ref{away-l} below.

We first recall a result providing useful bounds.

\llem{psnorm-l}{
There exist a constant $\varsigma>0$ such that,
for all $s \in \Re^m$ and all $v \geq 1$,
\beqn{s-sumsi}
\varsigma^v \|s\|^v
\leq \sum_{i\in \calN} \|s_i\|^v
\leq |\calN|  \,\|s\|^v.
\eeqn
}

\proof{See \cite[Lemma~4.1]{ChenToinWang17}.}

\noindent
This lemma states that $\sum_{i\in\calN}\|\cdot\|$ is a norm on $\Re^n$ whose
equivalence constants with respect to the Euclidean one are $\varsigma$ and
$|\calN|$.

Our next step is to specify under which conditions the standard
$\epsilon$-independent overestimation and derivative accuracy bounds typical of
the Lipschitz case (see \cite[Lemma~2.1]{CartGoulToin18b} for instance) can be
obtained for the elements functions of \req{problem}. We define, for a given
$k\geq 0$ and a given constant $\phi>0$ independent of $\epsilon$,
\beqn{Ok-def}
\calO_{k,\phi} \eqdef \{ i \in \calA_k^+
\mid \min[\, \|r_{i,k}\|, \,\|r_{i,k}+s_{i,k}\|\,] \geq \phi \}.
\eeqn
Observe that if, for some $i\in\calH$ and $b_i\not\in U_i\calF$, then
both $\|r_{i,k}\|$ and $\|r_{i,k}+s_{i,k}\|$ are bounded away from zero,
so $i \in \calO_{k,\phi}$ for all $k$ and all $\phi$ such that $\phi \leq
\min_{x\in\calF}\|U_ix-b_i\|$. Thus
we assume, without loss of generality, that
\beqn{binUF}
b_i\in U_i\calF \tim{ for all } i\in\calH.
\eeqn
We then obtain the following crucial error bounds.

\llem{Lip-th}{
Suppose that AS.2 and AS.3 hold. Then, for $k \geq 0$ and $L_{\max} \eqdef
\max_{i\in\calN}L_i$,
\beqn{fi-modi-N}
f_i(x_{i,k}+s_{i,k})
= m_i(x_{i,k},s_{i,k}) +
   \frac{1}{(p+1)!}\Big[ \tau_{i,k} L_{\max} -  \sigma_{i,k} \Big]\|s_{i,k}\|^{p+1}
\tim{ with }| \tau_{i,k} | \leq 1,
\eeqn
for all $i \in \calN$.
If, in addition,  $\phi >0$ is given and independent of
$\epsilon$, then there exists a constant $L(\phi)$ independent of $\epsilon$
such that, for $\ell \in \ii{p}$,
\beqn{gf-gmod}
\| \nabla_x^\ell f_{\calN\cup\calO_{k,\phi}}(x_k+s_k)
   -  \nabla_s^\ell T_{f_{\calN\cup\calO_{k,\phi}},p}(x_k,s_k) \|
\leq \frac{L(\phi)}{(p-\ell+1)!} \|s_k\|^{p-\ell+1}.
\eeqn
}

\proof{
  First note that, if $f_i$ has a Lipschitz continuous $p$-th derivative as a
  function of $x_i=U_ix$, then \req{taylor} shows that it also has a Lipschitz
  continuous $p$-th derivative as a function of $x$.  It is therefore enough to
  consider the element functions as functions of $x_i$.

  Observe now that, for each $k$ and $i \in \calN$, AS.2 and AS.3 ensure that
  \req{fi-modi-N} and the inequality
  \beqn{LipN-ders}
  \| \nabla_{x_i}^\ell f_i(x_{i,k}+s_{i,k})-\nabla_{s_i}^\ell T_{f_i,p}(x_{i,k},s_{i,k}) \|
  \leq \frac{L_i}{(p-\ell+1)!} \|s_{i,k}\|^{p-\ell+1}
  \eeqn
  immediately follow from the known bounds for $p$ times continuously
  differentiable functions with Lipschitz continuous $p$-th derivative (see
  \cite[Lemma~2.1]{CartGoulToin18b}). Consider now $i \in \calO_{k,\phi}$ for
  some $k$ and some fixed $\phi>0$, implying that
  $\min[\|r_{i,k}\|,\|r_{i,k}+s_{i,k}\|] \geq \phi >0$. Then
  \beqn{dersf+}
  \nabla_\cdot^\ell\|r_{i,k}+s_{i,k}\|^a[d]^\ell
  =\nabla_\cdot^\ell\big\| \|r_{i,k}+s_{i,k}\|u_{i,k}^+\big\|^a[d]^\ell
  = \nabla_\cdot^\ell\big\| \|r_{i,k}+s_{i,k}\|u_{i,k}\big\|^a[R_{i,k}d]^\ell
  \eeqn
  where $R_{i,k}$ is the rotation such that $R_{i,k}u_{i,k}^+=u_{i,k}$. We
  also have from \req{rotation} with $x$ replaced by $x_k+s_k$ that
  \beqn{dersT+}
  \nabla_{s_i}^\ell T_{m_i,p}(x_{i,k},s_{i,k})[d]^\ell
  = \nabla_\cdot^\ell \big\| \|r_{i,k}\|u_{i,k} \big\|^a[R_{i,k}d]^\ell.
  \eeqn
  Taking the difference between \req{dersf+} and \req{dersT+}, we obtain,
  successively using the definition of the tensor norm, the fact that
  $R_{i,k}$ is orthonormal and \req{difdersnorm} in Lemma~\ref{dersnorm-l},
  that
  \[
  \begin{array}{l}
  \big\|\nabla_\cdot^\ell\|r_{i,k}+s_{i,k}\|^a
  -\nabla_{s_i}^\ell T_{m_i,p}(x_{i,k},s_{i,k})\big\|_{[\ell]}\\*[1.5ex]
  \hspace*{30mm} = \bigmax_{\|d\|=1}\left|\nabla_\cdot^\ell\|r_{i,k}+s_{i,k}\|^a[d]^\ell
  -\nabla_{s_i}^\ell T_{m_i,p}(x_{i,k},s_{i,k})[d]^\ell\right|\\*[1.5ex]
  \hspace*{30mm} = \bigmax_{\|d\|=1}
  \left|\nabla_\cdot^\ell\big\| \|r_{i,k}+s_{i,k}\|u_{i,k}\big\|^a[R_{i,k}d]^\ell
      -\nabla_\cdot^\ell \big\| \|r_{i,k}\|u_{i,k} \big\|^a[R_{i,k}d]^\ell\right|\\*[1.5ex]
  \hspace*{30mm} =
  \left\|\,\nabla_\cdot^\ell\big\| \|r_{i,k}+s_{i,k}\|u_{i,k}\big\|^a
      -\nabla_\cdot^\ell \big\| \|r_{i,k}\|u_{i,k} \big\|^a\,\right\|_{[\ell]}\\*[1.5ex]
  \hspace*{30mm} =
  |\pi(a-\ell)| \left| \|r_{i,k}+s_{i,k}\|^{a-\ell}-\|r_{i,k}\|^{a-\ell}\right|.
  \end{array}
  \]
  Now the univariate function $\nu(t)\eqdef t^a$ is (more than) $p+1$ times
  continuously differentiable with bounded $(p+1)$-rst derivative on the
  interval $[t_1,t_2]$ and thus, from Lemma~\ref{Lip-th}, we have that
  \[
  \pi(a-\ell)\left|t_1^{a-\ell}-t_2^{a-\ell}\right|
  = \left|\frac{d^\ell\nu}{dt^\ell}(t_1)-\frac{d^\ell\nu}{dt^\ell}(t_2)\right|
  \leq \frac{L_\nu}{(p-\ell+1)!}|t_1-t_2|^{p-\ell+1},
  \]
  where $L_\nu$ is the upper bound on the $(p+1)$-rst derivative of $\nu(t)$ on
  interval $[t_1,t_2]$, that is
  $
  L_\nu = |\pi(a-p-1)|\min[t_1,t_2]^{a-p-1}.
  $
  As a consequence, we obtain that
  \[
  \big\|\nabla_\cdot^\ell\|r_{i,k}+s_{i,k}\|^a
  -\nabla_{s_i}^\ell T_{m_i,p}(x_{i,k},s_{i,k})\big\|_{[\ell]}
   \leq
  \bigfrac{L(\phi)}{(p-\ell+1)!}\big|\|r_{i,k}+s_{i,k}\|-\|r_{i,k}\|\big|^{p-\ell+1},
  \]
  where we use the fact that $\min[\|r_{i,k}\|,\|r_{i,k}+s_{i,k}\|] \geq \phi$
  to define
  \[
  L(\phi) = \max\big|\pi(a-p-1)|\phi^{a-p-1}, L_{\max}\big].
  \]
  We then observe that
  $\|s_{i,k}\|=\|r_{i,k}+s_{i,k}-r_{i,k}\|\geq
  \big|\|r_{i,k}+s_{i,k}\|-\|r_{i,k}\|\big|$
  which finally yields that
   \[
  \big\|\nabla_\cdot^\ell\|r_{i,k}+s_{i,k}\|^a
  -\nabla_{s_i}^\ell T_{m_i,p}(x_{i,k},s_{i,k})\big\|_{[\ell]}
  \leq \bigfrac{L(\phi)}{(p-\ell+1)!} \|s_{i,k}\|^{p-\ell+1}.
  \]
  Combining this last inequality with \req{LipN-ders} and the fact that
  $\nabla_{x_i}^\ell\|r_{i,k}+s_{i,k}\|^a=\nabla_\cdot^\ell\|r_{i,k}+s_{i,k}\|^a$
  then ensures that \req{gf-gmod} holds.
}

\noindent
Observe that the Lipschitz constant $L$ is independent of $\phi$ whenever
$\calH = \emptyset$. Our model definition also implies the following bound.

\llem{mdecr}{
  For all $k\geq 0$ before termination, $s_k \neq 0$, \req{rhok-def} is
  well-defined and
  \beqn{Dphi}
  \Delta T_k \geq \frac{\sigma_{\min}\varsigma^{p+1}}{(p+1)!}\, \|s_k\|^{p+1}.
  \eeqn
}

\proof{
  We immediately deduce that
  \beqn{Dphii}
  \Delta T_k \geq \frac{\sigma_{\min}}{(p+1)!}\, \sum_{i \in \calN}\|s_{i,k}\|^{p+1}
  \eeqn
  from \req{deltaT-def}, the observation that, at
  successful iterations, the algorithm enforces $\Delta m_k > 0$ and
  \req{sig-decr}. As a consequence, $s_k\neq 0$. Hence at least one
  $\|s_{i,k}\|$ is strictly positive because of \req{s-sumsi}, and \req{Dphii}
  therefore implies that \req{rhok-def} is well-defined. The inequality
  \req{Dphi} then follows from Lemma~\ref{psnorm-l}.
}

Following a now well-oiled track in convergence proofs for regularization
methods, we derive an upper bound on the regularization parameters.

\llem{sigmax}{
\cite[Lemma~4.6]{ChenToinWang17}
Suppose that AS.2 and AS.3 hold.  Then, for all $i \in \calN$ and all $k \geq 0$,
\beqn{sigma-max}
\sigma_{i,k} \in [\sigma_{\min}, \sigma_{\max}],
\eeqn
where $\sigma_{\max} \eqdef \gamma_2 L_{\max}$.
}

\proof{
Assume that, for some $i \in \calN$ and $k \geq 0$, $\sigma_{i,k} \geq L_i$.
Then \req{fi-modi-N} gives that \req{up-cond} must fail, ensuring
\req{sigma-max} because of the mechanism of the algorithm.
}

\noindent
It is important to note that $\sigma_{\max}$ is independent of
$\epsilon$. We now verify that the trial step produced by Step~2 of the
PSAR$p$ Algorithm either essentially fixes the residuals $r_i$ to zero (their
value being then fixed for the rest of the calculation), or is long enough
(i.e.\ \req{longs} holds), or maintains these residuals safely away from zero
in the sense that their norm exceeds an $\epsilon$-independent constant.

\llem{away-l}{
Suppose that AS.1, AS.2, AS.3 and AS.5 hold, that $\calH \neq
\emptyset$ and that \req{longs} fails. Let
\beqn{omega-def}
\omega \eqdef
\min\left[ \frac{a}{16},
  \left(\frac{a}{12|\calN|\left(\kappa_\calN
    +\frac{\sigma_{\max}}{(p-q+1)!}\right)}\right)^{\frac{1}{1-a}}\right].
\eeqn
Then, if, for some $i \in \calH$,
\beqn{a-assxik}
\|r_{i,k}\| < \omega,
\eeqn
we have that
\beqn{gap}
\|r_{i,k}+ s_{i,k}\| \leq \epsilon
\tim{ or }
\|r_{i,k}+ s_{i,k}\| \geq \omega.
\eeqn
}

\proof{The conclusion is obvious if $i \in \calC_k^+= \calH\setminus\calA_k^+$.
Consider now $i \in \calA_k^+$ and suppose, for the purpose of deriving a
contradiction, that
\beqn{absurd}
\|r_{i,k}+s_{i,k}\| \in (\epsilon, \omega) \tim{ for some } i \in \calA_k^+,
\eeqn
and immediately note that
the failure of \req{longs} and the orthonormality of the rows of $U_i$ imply that
\beqn{skleq1}
\|s_{i,k}\| \le \|s_k\|< \varpi \epsilon^{\frac{1}{p-q+1}} \leq 1
\eeqn
and also that \req{term-m} must hold.
As a consequence, for some $\delta_k \in (0,1]$,
\beqn{psib-1a}
a \|r_{i,k}+s_{i,k}\| \chi_q(\delta_k)
\geq \psi_{m,q}^{\epsilon,\delta_k}(x_k,s_k).
\eeqn
Consider now the vector
\beqn{dk-def}
d_k = -\min\big[\delta_k ,\|r_{i,k}+s_{i,k}\| \big] v_{i,k}^+
\tim{with} v_{i,k}^+ = U_i^\dagger u_{i,k}^+\eqdef U_i^\dagger
\frac{r_{i,k}+s_{i,k}}{\|r_{i,k}+s_{i,k}\|}.
\eeqn
We now verify that $d_k$ is admissible for problem \req{psim-def}. Clearly
$\|d_k\| = \delta_k$ because the rows of $U_i$ orthonormal.
We also see that \req{Ui-ortho-H} and \req{Rii-def} imply that, since $i \in
\calA_k^+$,
\beqn{dkinRk+}
d_k \in \calR_{\{i\}} \subseteq \calR_k^+.
\eeqn
Moreover, we have that
\beqn{ininter}
x_k+s_k+d_k \in [\![ x_k+s_k, x_k+s_k- U_i^\dagger(r_{i,k}+s_{i,k})]\!],
\eeqn
where $[\![v,w]\!]$ denotes the line segment joining the vectors $v$ and $w$.
But
\[
\begin{array}{ll}
x_k+s_k - U_i^\dagger(r_{i,k}+s_{i,k})
& = x_k+s_k - U_i^\dagger U_i(x_k+s_k) + U_i^\dagger b_i \\
& = (I- U_i^\dagger U_i)(x_k+s_k) + U_i^\dagger b_i \\
& = P_{\ker(U_i)} [x_k+s_k] + U_i^\dagger b_i\\
& \in \calF,
\end{array}
\]
where we have used \req{kernel-centered} to deduce the last inclusion.
Since $\calF$ is convex and $x_k+s_k\in \calF$, we deduce from
\req{ininter} that $x_k+s_k +d_k \in \calF$. As a consequence, $d_k$ is
admissible for problem \req{psim-def} and hence, using \req{psib-1a},
\beqn{psib-1}
a \|r_{i,k}+s_{i,k}\| \chi_q(\delta_k)
\geq \psi_{m,q}^{\epsilon,\delta_k}(x_k,s_k)
\geq \max\left[0,m(x_k,s_k)-T_{m,q}(x_k,s_k-d_k)\right].
\eeqn
Moreover \req{dkinRk+} and \req{miH-def} imply that
\beqn{psib-2}
\begin{array}{ll}
m(x_k,s_k)-T_{m,q}(x_k,s_k-d_k)&\\*[1.5ex]
&\hspace*{-4cm} =
   m_\calN(x_k,s_k)-T_{m_\calN,q}(x_k,s_k-d_k)
   + m_i(x_{i,k},s_{i,k})- m_i(x_{i,k},s_{i,k}-U_id_k)\\*[1.5ex]
&\hspace*{-4cm} \geq
   -\left|m_\calN(x_k,s_k)-T_{m_\calN,q}(x_k,s_k-d_k)\right|
   + m_i(x_{i,k},s_{i,k})- m_i(x_{i,k},s_{i,k}-U_id_k).
\end{array}
\eeqn
We start by considering the first term in the right-hand side of this
inequality. Observe now that \req{a-assxik} ensures that $x_k\in \calV$ (as
defined in \req{calV-def}). Hence AS.5, \req{a-assxik} and \req{skleq1}
together imply that, for each $i\in \calN$,
\beqn{psib-3}
\begin{array}{l}
|m_i(x_k,s_k)-T_{m_i,q}(x_k,s_k-d_k)|\\*[2ex]
\hspace*{10mm}
\leq \left|\bigsum_{\ell=1}^q \frac{1}{\ell!} \nabla_x^\ell
T_{m_i,q}(x_{i,k},s_{i,k})[-U_id_k]^\ell\right|\\*[2.5ex]
\hspace*{10mm}
= \left|\bigsum_{\ell=1}^q \frac{1}{\ell!}\left(
\sum_{t=\ell}^p \frac{1}{(t-\ell)!}\nabla_x^t f_i(x_{i,k})[s_{i,k}]^{t-\ell}
   + \bigfrac{\sigma_{i,k}}{(p+1)!}\big\|\nabla_\cdot^\ell\|s_{i,k}\|^{p+1}
   \big\| \right)[-U_id_k]^\ell\right|.
\end{array}
\eeqn
Using now the identity $\|U_id_k\| = \|d_k\|= \delta_k$ and the fact that
\[
\sum_{t=\ell}^p \frac{1}{(t-\ell)!} \leq 1+\chi_{p-\ell}(1) < 3,
\]
we obtain from \req{psib-3}, the triangle inequality and \req{skleq1} that
\beqn{psib-4}
|m_i(x_k,s_k)-T_{m_i,q}(x_k,s_k-d_k)|
< \bigsum_{\ell=1}^q\frac{1}{\ell!} \left(3 \|\nabla_x^t f_i(x_{i,k})\|_{[t]}
   + \bigfrac{\sigma_{i,k}}{(p+1)!}\big\|\nabla_\cdot^\ell\|s_{i,k}\|^{p+1}\big\|\right)\delta_k^\ell.
\eeqn
But we have from Lemma~\ref{dersnorm-l} and \req{skleq1} that, for
$\ell\in \ii{q}$,
\beqn{psib-5}
\big\|\nabla_\cdot^\ell\|s_{i,k}\|^{p+1}\big\|
= |\pi(p-\ell+1)|\,\|s_{i,k}\|^{p+1-\ell}
\leq \frac{(p+1)!}{(p-q+1)!},
\eeqn
and therefore that
\beqn{psib-6}
\begin{array}{ll}
\left|m_\calN(x_k,s_k)-T_{m_\calN,q}(x_k,s_k-\delta_k \|r_{i,k}+s_{i,k}\|v_{i,k}^+)\right|
& < 3|\calN|\left(\kappa_\calN+\frac{1}{(p-q+1)!}\sigma_{\max}\right)\chi_q(\delta_k)\\*[1.5ex]
& \leq \quarter a \omega^{a-1}\chi_q(\delta_k),
\end{array}
\eeqn
where we have used \req{omega-def} to derive the last inequality.
Let us now consider the second term in the right-hand side of \req{psib-2}.
Applying Lemma~\ref{miprop-l}, we obtain that
$\mu(\|r_{i,k}+s_{i,k}\|,\cdot)$ is concave between $0$ and $-\|r_{i,k}+s_{i,k}\|$
and $\mu(\|r_{i,k}\|,\cdot)$  is concave between $0$ and $-\|r_{i,k}\|$.
Therefore, because \req{a-assxik}, we may deduce that
\[
\begin{array}{ll}
m_i(x_{i,k},s_{i,k})- m_i(x_{i,k},s_{i,k}-U_id_k)
& = \mu(\|r_{i,k}+s_{i,k}\|,0) - \mu(\|r_{i,k}+s_{i,k}\|,\|U_id_k\|)\\*[1.5ex]
&\geq \nabla_\zeta^1\mu(\|r_{i,k}+s_{i,k}\|,0)\|U_id_k\|\\*[1.5ex]
&\geq \nabla_\zeta^1\mu(\|r_{i,k}\|,\|r_{i,k}+s_{i,k}\|-\|r_{i,k}\|)\|U_id_k\|\\*[1.5ex]
&\geq a\|r_{i,k}\|^{a-1}\delta_k\\*[1.5ex]
&\geq \half a \omega^{a-1}\chi_q(\delta_k),
\end{array}
\]
where the second and third inequalities result from \req{dersmi-odd}.
Combining now this inequality with \req{psib-1}, \req{psib-2} and \req{psib-6}, we deduce
that
\[
a\|r_{i,k}+s_{i,k}\|\chi_q(\delta_k)
> \half a \omega^{a-1}\chi_q(\delta_k) - \quarter a \omega^{a-1}\chi_q(\delta_k)
= \quarter a \omega^{a-1} \chi_q(\delta_k).
\]
Finally, we obtain using \req{absurd} that
\[
\omega > \quarter \omega^{a-1},
\]
which impossible in view of \req{omega-def}.  Hence \req{absurd} cannot hold
and the proof is complete.
}

\noindent
This last result is crucial in that it shows that there is a ``forbidden''
interval $(\epsilon,\omega)$ for the residual's norms $\|r_i(x_k+s_k)\|$,
where $\omega$ only depends on the problem and is independent of $\epsilon$.  This in turn
allows to partition the successful iterates into subsets, distinguishing
iterates which ``fix'' a residual to a near zero value, iterates with long
steps and iterates with possibly short steps in regions where the considered
objective function's $p$-th derivative tensor is safely bounded independently
of $\epsilon$.  Our analysis now follows the broad outline of
\cite{ChenToinWang17} while simplifying some arguments. Focusing on the case
where $\calH \neq \emptyset$, we first isolate the set of successful iterations which
``deactivate'' a residual, that is
\[
\calS_\epsilon \eqdef \{ k \in \calS \mid
                  \|r_{i,k}+s_{i,k}\| \leq \epsilon
                  \tim{and}
                  \|r_{i,k}\| > \epsilon
                  \tim{for some} i \in \calH\},
\]
and notice that, by construction
\beqn{Seps-bound}
|\calS_\epsilon| \leq |\calH|.
\eeqn
We next define the $\epsilon$-independent constant
\[
\alpha = \threequarters \omega
\]
and
\beqn{Ss-def}
\calS_{\|s\|}
\eqdef \{ k \in \calS \mid \|s_k\| \geq \quarter \omega \}.
\eeqn
Moreover, for an iteration $k\in \calS\setminus (\calS_\epsilon \cup
\calS_{\|s\|})$, we verify that $\calA_k$ can be partitioned into
\[
\begin{array}{ll}
\calI_{\heartsuit,k} \eqdef
& \{ i \in \calA_k \mid \|r_{i,k}\| \in [\alpha,+\infty) \tim{and}
                        \|r_{i,k}+s_{i,k}\| \in [\alpha,+\infty)\} \\
\calI_{\diamondsuit,k} \eqdef
& \{ i \in \calA_k \mid \big( \|r_{i,k}\| \in [\omega,+\infty) \tim{and}
                             \|r_{i,k}+s_{i,k}\| \in (\epsilon,\alpha) \big)\\
& \hspace*{13mm}        \tim{ or }
                        \big( \|r_{i,k}\| \in (\epsilon,\alpha) \tim{and}
                             \|r_{i,k}+s_{i,k}\| \in [\omega,\infty) \big)
                               \} \\
\calI_{\clubsuit,k} \eqdef
& \{ i \in \calA_k \mid \|r_{i,k}\| \in (\epsilon,\omega) \tim{and}
                        \|r_{i,k}+s_{i,k}\| \in (\epsilon,\omega)\}.
\end{array}
\]
Morever, Lemma~\ref{away-l} shows that $\calI_{\clubsuit,k}$ is always empty
and one additionally has that, if $i\in\calI_{\diamondsuit,k}$, then
\[
\|s_k\|
\geq \|s_{i,k}\|
\geq \big|\|r_{i,k}+s_{i,k}\|-\|r_{i,k}\|\big|
\geq \omega - \alpha
= \quarter \omega,
\]
implies that $k \in \calS_{\|s\|}$.  Hence $\calI_{\diamondsuit,k}$ is also
empty and
\beqn{Bk-heart}
\calA_k = \calI_{\heartsuit,k}
\tim{ for } k\in \calS\setminus (\calS_\epsilon\cup\calS_{\|s\|}) \eqdef \calS_\heartsuit.
\eeqn

The next important result shows that steps at iteration belonging to
$\calS_\heartsuit$ are long enough, because they are taken over region where a
good $\epsilon$-independent Lipschitz bounds holds.  Indeed, if
$\calH\neq\emptyset$ and assuming that $\epsilon\leq \alpha$, we have, for $k
\in \calS_\heartsuit$, that $\calA_k^+ = \calA_k$ and thus that
$\calW_k^+=\calW_k$ and $\calR_k^+=\calR_k$. Moreover, the definition of
$I_{\heartsuit,k}  = \calA_k$ ensures that
$
\calA_k \subseteq \calO_{k,\alpha}
$
and thus that Lemma~\ref{Lip-th} (and in particular \req{gf-gmod}) guarantees
that $f_{\calW_k}$ satisfies standard derivative error bounds for functions
with Lipschitz continuous $p$-th derivative (with corresponding Lipschitz
constant $L(\alpha)$). We may therefore apply known results for such functions
to $f_{\calW_k}$.  The following lemma is extracted from \cite{CartGoulToin18b}, by
specializing Lemma~3.3 in that reference to optimization of $f_{\calW_k}$ over
$\calR_k$ for functions with Lipschitz continuous $p$-th derivative
(i.e. $\beta=1$ in \cite{CartGoulToin18b}).

\llem{longs-Lip-l}{
  Suppose that AS.1 -- AS.3 and AS.5 hold, that
  \beqn{epsalpha}
  \epsilon \leq \alpha \tim{ if } \calH \neq \emptyset
  \eeqn
  and consider $k \in \calS_\heartsuit$ such that the PSAR$p$ Algorithm does
  not terminate at iteration $k+1$. Then
  \beqn{longs-Lip}
  \|s_k\| \geq \kappa_\heartsuit \epsilon^{\frac{1}{p-q+1}}
  \tim{ with }
  \kappa_\heartsuit
  \eqdef \left(\frac{(p-q+1)!}{L(\alpha)+\theta +\sigma_{\max}}\right)^{\frac{1}{p-q+1}}.
  \eeqn
}

\noindent
We may finally establish our final evaluation complexity bound by combining
our results so far.

\lthm{compl-succ-l}{Suppose that AS.1--AS.5 and \req{epsalpha} hold.
  Then the PSAR$p$ Algorithm requires at most
  \beqn{max-succ-its}
  \left\lfloor \kappa_\calS(f(x_0)-f_{\rm low}) \epsilon^{-\frac{p+1}{p-q+1}}\right\rfloor+ |\calH|
  \eeqn
  successful iterations and at most
  \beqn{max-evals}
  \left\lfloor \left \lfloor
   \kappa_S ( f(x_0)- f_{\rm low})
   \left( \epsilon^{-\frac{p+1}{p-q+1}} \right)
   \right \rfloor
                 \left(1+\frac{|\log\gamma_1|}{\log\gamma_2}\right)+
  \frac{1}{\log\gamma_2}\log\left(\frac{\sigma_{\max}}{\sigma_0}\right)\right\rfloor
  + |\calH| +1
  \eeqn
  evaluations of $f$ and its $p$ first derivatives to return an
  $(\epsilon,\delta)$-approximate $q$-th-order-necessary minimizer for problem
  \req{problem}, where
  \beqn{kappaS-def}
  \kappa_S \eqdef \frac{(p+1)!}{\eta\sigma_{\min}\varsigma^{p+1}}
        \left(\frac{(p-q+1)!}
                   {L(\alpha)+\theta +\sigma_{\max}}\right)^{-\frac{1}{p-q+1}}.
  \eeqn
}
\proof{
  Consider $k\in \calS$ before termination.
  Because the iteration is successful, we obtain from AS.4, Step~3
  of the algorithm and Lemma~\ref{mdecr} that
  \beqn{cs-1}
  f(x_0)-f_{\rm low}
  \geq f(x_0)-f(x_{k+1})
  = \sum_{j\in \calS_k}\Delta f_k
  \geq \eta \sum_{j\in \calS_k} \Delta T_k
  \geq \frac{\eta\sigma_{\min}\varsigma^{p+1}}{(p+1)!}
       \sum_{j\in \calS_k}\|s_k\|^{p+1}.
  \eeqn
  Defining now
  \[
  \calS_{\epsilon,k} \eqdef \calS_\epsilon \cap \iibe{0}{k},
  \ms
  \calS_{\|s\|,k} \eqdef \calS_{\|s\|} \cap \iibe{0}{k}
  \tim{and}
  \calS_{\heartsuit,k} \eqdef \calS_\heartsuit \cap \iibe{0}{k},
  \]
  we verify that $\calS_{\|s\|,k}$ and $\calS_{\heartsuit,k}$ form a partition
  of $\calS_k\setminus \calS_{\epsilon,k}$.
  As a consequence, we have that
  \[
  \begin{array}{ll}
      f(x_0)-f_{\rm low}
      & \geq \bigfrac{\eta\sigma_{\min}\varsigma^{p+1}}{(p+1)!}
      \left\{
      |\calS_{\|s\|,k}|\bigmin_{j \in \calS_{\|s\|,k}}\|s_k\|^{p+1}
      + |\calS_{\heartsuit,k}|\bigmin_{j \in \calS_{\heartsuit,k}}\|s_k\|^{p+1}
      \right\}\\
      & \geq \bigfrac{\eta\sigma_{\min}\varsigma^{p+1}}{(p+1)!}
      \left\{
      |\calS_{\|s\|,k}|(\quarter \omega)^{p+1}
      + |\calS_{\heartsuit,k}|\left(\kappa_\heartsuit \epsilon^{\frac{1}{p-q+1}}\right)^{p+1}
      \right\}\\
      & \geq \bigfrac{\eta\sigma_{\min}\varsigma^{p+1}}{(p+1)!}
      \left\{|\calS_{\|s\|,k}|+|\calS_{\heartsuit,k}|\right\}
      \min\left[(\quarter \omega)^{p+1},\left(\kappa_\heartsuit \epsilon^{\frac{1}{p-q+1}}\right)^{p+1}\right]\\
      & \geq \bigfrac{\eta\sigma_{\min}\varsigma^{p+1}}{(p+1)!}
      |\calS_k\setminus\calS_{\epsilon,k}|\,\kappa_\heartsuit^{\frac{1}{p-q+1}}\,\epsilon^{\frac{1}{p-q+1}},
  \end{array}
  \]
  where we have used \req{cs-1}, Lemma~\ref{psnorm-l}, \req{Ss-def} and
  \req{longs-Lip} to deduce the second inequality, and the assumption that
  (without loss of generality in view of \req{longs-Lip})
  $\kappa_\heartsuit\leq \quarter \omega$ to deduce the last. The above
  inequality yields that
  \[
  |\calS_k|
  = |\calS_k\setminus\calS_{\epsilon,k}| + |\calS_{\epsilon,k}|
  \leq \kappa_S (f(x_0)-f_{\rm low})\epsilon^{-\frac{p+1}{p-q+1}}+|\calS_{\epsilon,k}|,
  \]
  where $\kappa_S$ is given by \req{kappaS-def}.
  Since $|\calS_{\epsilon,k}|\leq |\calS_\epsilon| \leq |\calH|$, we finally
  deduce that the bound \req{max-succ-its} holds.  The bound \req{max-evals}
  then follows by applying Lemma~\ref{a-succ-unsucc} and observing that $f$ and its first
  $p$ derivatives are evaluated at most once per iteration, plus once at termination.
}

We conclude our development by recalling that the above result is valid for
$\calH = \emptyset$, in which case the problem is a smooth
convexly-constrained partially-separable problem. Note that the
norm-equivalence constant $\varsigma$ occurs in \req{kappaS-def}, which
indicate that the underlying geometry of the problem's invariants subspaces
$\ker(U_i)$ may have a significant impact on complexity.

\numsection{Conclusions}\label{concl-s}

We have shown that an $(\epsilon,\delta)$-approximate
$q$-th-order critical point of partially-separable convexly-constrained
optimization with non-Lipschitzian singularities can be found at most
$O(\epsilon^{-(p+1)/(p-q+1)})$ evaluations of the objective function and its
first $p$ derivatives for any $q \in \{1,2,\ldots, p\}$ whenever the smooth
element functions $f_i$, $i\in \calN$ of the objective function are $p$ times
differentiable. This worst-case complexity is obtained via our Algorithm 4.1
(PSAR$_p$) with an $p$-th order Taylor model which uses the underlying
rotational symmetry of the Euclidean norm function for $f_\calH$ and the first
$p$ derivatives (whenever they exist) of the ``element functions'' $f_i$, for
$i\in \calM$.

Several observations are of interest. A first one is that the results remain
valid if Lipschitz continuity is not assumed on the whole of the feasible set,
but restricted to the segments of the ``path of iterates'', that is
$\cup_k[\![x_k,x_{k+1}]\!]$. While this might in general be difficult to ensure a
priori, there may be case where problem structure could help. A second
observation is that convexity of the feasible set is only used on the segments
$\cup_{i,k} [\![x_{i,k}, U_i^\dagger b_i]\!]$.  Again this might be exploitable in
some cases. The third observation is that, in line with \cite{CartGoulToin18b},
it is possible to replace the Lipschitz continuity assumption by a weaker
H\"{o}lder continuity.

While it may be possible to handle non-kernel-centered
feasible sets (maybe along the lines of the discussion in
\cite{ChenToinWang17}), this remains open at this stage.  Another interesting
perspective is a more generic exploitation of geometric symmetries inherent to
optimization problems: our treatment here focuses on a specific case of
rotational symmetry, but this should not, one hopes, be limitative.

\noindent
\section*{\footnotesize Acknowledgements}
\vspace*{-3mm}
{\footnotesize
Xiaojun Chen would like to thank Hong Kong Research Grant Council for grant
PolyU153000/17p. Philippe Toint would like to thank the Hong Kong Polytechnic
University for its support while this research was being conducted.


\begin{thebibliography}{9}
\bibitem{Ahsen}
M.E. Ahsen and M. Vidyasagar, Error bounds for compressed sensing algorithms with group sparsity: A unified approach,
Appl. Comput. Harmon. Anal. 43(2017), 212-232.

\bibitem{beck2017}
A. Beck and N. Hallak, Optimization problems involving group sparsity terms,
  Math. Program., (2018), online.


\bibitem{balda}
L. Baldassarre, N. Bhan, V. Cevher, A. Kyrillidis and S. Satpathi, Group-sparse model selection: hardness and relaxations, IEEE Trans. Inf. Theory, 62(2016), 6508-6534.

\bibitem{BellGuriMoriToin18}
S.~Bellavia, G.~Gurioli, B.~Morini, and {Ph.}~L. Toint.
\newblock Deterministic and stochastic inexact regularization algorithms for
  nonconvex optimization with optimal complexity.
\newblock arXiv:1811.03831, 2018.

\bibitem{Bian_Chen_SIOPT}
W. Bian and  X. Chen, Worst-case complexity of smoothing quadratic regularization methods for non-Lipschitzian optimization, SIAM J. Optim., 23(2013), 1718-1741.

\bibitem{Bian-Chen-Ye}
W. Bian, X. Chen and Y. Ye, Complexity analysis of interior point algorithms for non-Lipschitz and nonconvex minimization,  Math. Program., 149(2015), 301-327.

\bibitem{BirgGardMartSantToin17}
E.~G. Birgin, J.~L. Gardenghi, J.~M. Mart\'{i}nez, S.~A. Santos, and Ph.~L.
  Toint.
\newblock Worst-case evaluation complexity for unconstrained nonlinear
  optimization using high-order regularized models.
\newblock Math. Program., 163(1):359--368, 2017.

\bibitem{breheny2015}
P. Breheny and J. Huang, Group descent algorithms for nonconvex penalized
linear and logistic regression models with grouped predictors,
Stat. Comput., 25(2015), 173-187.

\bibitem{CartGoulToin12b}
C.~Cartis, N.~I.~M. Gould, and Ph.~L. Toint.
\newblock An adaptive cubic regularization algorithm for nonconvex optimization
  with convex constraints and its function-evaluation complexity.
\newblock J. Numer. Anal., 32(4):1662--1695, 2012.

\bibitem{CartGoulToin17c}
C.~Cartis, N.~I.~M. Gould, and Ph.~L. Toint.
\newblock Second-order optimality and beyond: characterization and evaluation
  complexity in convexly-constrained nonlinear optimization.
\newblock Found. Comp. Math., 18(5):1073--1107,
  2018.

\bibitem{CartGoulToin18b}
C.~Cartis, N.~I.~M. Gould, and Ph.~L. Toint.
\newblock Sharp worst-case evaluation complexity bounds for arbitrary-order
  nonconvex optimization with inexpensive constraints.
\newblock arXiv:1811.01220, 2018.

\bibitem{CartGoulToin18a}
C.~Cartis, N.~I.~M. Gould, and Ph.~L. Toint.
\newblock Worst-case evaluation complexity and optimality of second-order
  methods for nonconvex smooth optimization.
\newblock To appear in the Proceedings of the 2018 International Conference of
  Mathematicians (ICM 2018), Rio de Janeiro, 2018.

\bibitem{ChenDengZhan98}
L.~Chen, N.~Deng, and J.~Zhang.
\newblock Modified partial-update {N}ewton-type algorithms for unary
  optimization.
\newblock  J. Optim. Theory  Appl.,
  97(2):385--406, 1998.

\bibitem{CXY}
X. Chen, F. Xu and Y. Ye, Lower bound theory of nonzero
entries in solutions of $\ell_2$-$\ell_p$ minimization,  SIAM J.
Sci. Comput., 32(2010), 2832-2852.


\bibitem{CNY}
 X. Chen, L. Niu and Y. Yuan,  Optimality conditions and smoothing trust region Newton method for non-Lipschitz optimization,
SIAM J. Optim., 23(2013), 1528-1552.

\bibitem{CGWY}
X. Chen, D. Ge, Z. Wang and Y. Ye, Complexity of unconstrained
$L_2$-$L_p$ minimization,  Math. Program., 143(2014), 371-383.

\bibitem{ChenToinWang17}
X.~Chen, Ph.~L. Toint, and H.~Wang.
\newblock  Complexity of partially-separable convexly-constrained optimization with non-{L}ipschitzian singularities,
\newblock to appear in  SIAM J. Optim.

\bibitem{Chen_Rob}
  X. Chen  and R. Womersley, Spherical designs and nonconvex minimization for recovery of sparse signals on the sphere,  SIAM J. Imaging Sci., 11(2018), 1390-1415.


\bibitem{ConnGoulSartToin96a}
A.~R. Conn, N.~I.~M. Gould, A.~Sartenaer, and Ph.~L. Toint.
\newblock Convergence properties of minimization algorithms for convex
  constraints using a structured trust region.
\newblock  SIAM J. Optim., 6(4):1059--1086, 1996.

\bibitem{ConnGoulToin92}
A.~R. Conn, N.~I.~M. Gould, and Ph.~L. Toint.
\newblock {\em {\sf LANCELOT}: a {F}ortran package for large-scale nonlinear
  optimization ({R}elease {A})}.
\newblock Number~17 in Springer Series in Computational Mathematics. Springer
  Verlag, Heidelberg, Berlin, New York, 1992.

\bibitem{ConnGoulToin00}
A.~R. Conn, N.~I.~M. Gould, and Ph.~L. Toint.
\newblock {\em Trust-Region Methods}.
\newblock MPS-SIAM Series on Optimization. SIAM, Philadelphia, USA, 2000.

\bibitem{Eldar}
Y.C. Eldar, P. Kuppinger and H. B\"olcskei, Block-sparse signals: uncertainty relations and efficient recovery,
 IEEE Trans. Signal Process,
  58(2010), 3042-3054.

\bibitem{FourGayKern87}
R.~Fourer, D.~M. Gay, and B.~W. Kernighan.
\newblock {AMPL}: A mathematical programming language.
\newblock Computer science technical report, AT\&T Bell Laboratories, Murray
  Hill, USA, 1987.

\bibitem{Gay96}
D.~M. Gay.
\newblock Automatically finding and exploiting partially separable structure in
  nonlinear programming problems.
\newblock Technical report, Bell Laboratories, Murray Hill, New Jersey, USA,
  1996.

\bibitem{GoldWang93}
D.~Goldfarb and S.~Wang.
\newblock Partial-update {N}ewton methods for unary, factorable and partially
  separable optimization.
\newblock SIAM J.  Optim., 3(2):383--397, 1993.

\bibitem{GoulOrbaToin15b}
N.~I.~M. Gould, D.~Orban, and Ph.~L. Toint.
\newblock {\sf CUTEst}: a constrained and unconstrained testing environment
  with safe threads for mathematical optimization.
\newblock Comp. Optim. Appl., 60(3):545--557,
  2015.

\bibitem{GoulToin07b}
N.~I.~M. Gould and Ph.~L. Toint.
\newblock {\sf {FILTRANE}}, a {F}ortran~95 filter-trust-region package for
  solving systems of nonlinear equalities, nonlinear inequalities and nonlinear
  least-squares problems.
\newblock ACM Trans. Math. Soft., 33(1):3--25, 2007.

\bibitem{GrieToin82a}
A.~Griewank and Ph.~L. Toint.
\newblock On the unconstrained optimization of partially separable functions.
\newblock In M.~J.~D. Powell, editor, {\em Nonlinear Optimization 1981}, pages
  301--312, London, 1982. Academic Press.

\bibitem{huang2009}
J. Huang, S.  Ma, H. Xie and C. Zhang, A group bridge approach for variable selection,
  Biometrika, 96(2009), 339-355.

\bibitem{huang2010benefit}
J. Huang and  T. Zhang,  The benefit of group sparsity,
  Ann. Stat., 38(2010), 1978-2004.
\bibitem{gia2018}
G. Le, I. Sloan, R. Womersley and Y. Wang, Sparse isotropic regularization for spherical Harmonic representations of random fields on the sphere,
  arXiv preprint arXiv:1801.03212, (2018).

    \bibitem{juditsky2012}
A. Juditsky, F. Karzan, A. Nemirovski and B. Polyak, Accuracy guaranties for $\ell_1 $ recovery of block-sparse signals,  Ann. Stat., 40(2012), 3077-3107.

  \bibitem{Lee2012}
  K. Lee, Y. Bresler and M. Junge, Subspace methods for joint sparse recovery, IEEE Trans. Inform. Theory, 58(2012),
  3613-3641.

  \bibitem{Lee2016}
S. Lee, M. Oh and Y. Kim, Sparse optimization for nonconvex group penalized estimation,
  J. Stat. Comput. Simulat, 86(2016), 597-610.

  \bibitem{Lv}
  X. Lv, G. Bi and C. Wan, The group Lasso for stable recovery of block-sparse signal representations,
  IEEE Trans. Signal Proc., 59(2011), 1371-1382.

\bibitem{MaHuang}
S. Ma and J. Huang, A concave pairwise fusion approach to subgroup
analysis, J. Amer. Stat. Assoc., 112(2017), 410-423.


\bibitem{MareRichTaka14}
J.~Mare\v{c}ek, P.~Richt\'arik, and M.~Tak\'a\v{c}.
\newblock Distributed block coordinate descent for minimizing partially
  separable functions.
\newblock Technical report, Department of Mathematics and Statistics,
  University of Edinburgh, Edinburgh, Scotland, 2014.
\bibitem{obozinski2011}
G. Obozinski, M.J. Wainwright and  M. Jordan, Support union recovery in high-dimensional multivariate regression, Ann. Statist., 39(2011), 1-47.

 \bibitem{Yuan}
  M. Yuan and Y. Lin, Model selection and estimation in regression with grouped variables, J. R. Statist. Soc. B 68(2006), 49-67.

\end{thebibliography}
}

\appendix

\appnumsection{Appendix}

\noindent
\textbf{Proof Lemma~\ref{dersnorm-l}}
The proof of \req{dersnorm} is essentially borrowed from
\cite[Lemma~2.4]{CartGoulToin18b}, although details differ because the
present version covers $a \in (0,1)$. We first observe that $\nabla_\cdot^j
\|r\|^a$ is a $j$-th order tensor,
whose norm is defined using \req{Tnorm}.  Moreover, using the relationships
\beqn{donce}
\nabla_\cdot^1 \|r\|^\tau  = \tau\, \|r\|^{\tau-2}r
\tim{ and }
\nabla_\cdot^1 \big(r^{\tau \otimes}\big) = \tau\, r^{(\tau-1)\otimes}\otimes I,
\ms (\tau \in \Re),
\eeqn
defining
\beqn{munu-def}
\nu_0 \eqdef 1,
\tim{ and }
\nu_i \eqdef \prod_{\ell=1}^{i}(a+2-2\ell),
\eeqn
and proceeding by induction, we obtain that, for some $\mu_{j,i}\geq 0$ with
$\mu_{1,1}=1$,
\[
\begin{array}{ll}
\nabla_\cdot^1\left[\nabla_\cdot^{j-1} \| r\|^a  \right]&\\
& \hspace*{-1cm} = \nabla_\cdot^1\left[ \bigsum_{i=2}^j \mu_{j-1,i-1} \nu_{i-1}
   \|r\|^{a-2(i-1)} \, r^{(2(i-1)-(j-1)) \otimes} \otimes I^{((j-1)-(i-1))\otimes} \right]\\
&\hspace*{-1cm} = \bigsum_{i=2}^j \mu_{j-1,i-1} \nu_{i-1} \Big[
  (a-2(i-1))\|r\|^{a-2(i-1)-2} \, r^{(2(i-1)-(j-1)+1) \otimes} \otimes I^{(j-i)\otimes}\\
&  + ((2(i-1)-(j-1)) \|r\|^{a-2(i-1)} \, r^{(2(i-1)-(j-1)-1)\otimes}
  \otimes I^{(j-1)-(i-1)+1)\otimes} \Big]\\
&\hspace*{-1cm} = \bigsum_{i=2}^j \mu_{j-1,i-1} \nu_{i-1} \Big[
  (a+2-2i)\|r\|^{a-2i} \, r^{(2i-j) \otimes} \otimes I^{(j-i)\otimes}\\
  &  + (2(i-1)-j+1) \|r\|^{a-2(i-1)} \, r^{(2(i-1)-j) \otimes}
  \otimes I^{(j-(i-1))\otimes} \Big]\\
&\hspace*{-1cm} = \bigsum_{i=2}^j \mu_{j-1,i-1} \nu_{i-1}
  (a+2-2i)\|r\|^{a-2i} \, r^{(2i-j) \otimes} \otimes I^{(j-i)\otimes}\\
&  + \bigsum_{i=1}^{j-1} (2i-j+1) \mu_{j-1,i}\nu_i \|r\|^{a-2i}\, r^{(2i-j) \otimes}
  \otimes I^{(j-i)\otimes} \\
& \hspace*{-1cm} = \bigsum_{i=1}^j\big((a+2-2i)\mu_{j-1,i-1}\nu_{i-1}
   + (2i-j+1)\mu_{j-1,i}\nu_i \big) \|r\|^{a-2i} \, r^{(2i-j) \otimes} \otimes I^{(j-i)\otimes},
\end{array}
\]
where the last equation uses the convention that $\mu_{j,0} = 0$ and $\mu_{j-1,j} = 0$ for all $j$.
Thus we may write
\beqn{nablaj-ns-full}
\nabla_\cdot^j \|r\|^a
=\nabla_\cdot^1\left[\nabla_\cdot^{j-1}\| r\|^a \right]
= \bigsum_{i=1}^j \mu_{j,i} \nu_i \,
   \|r\|^{a-2i} \, r^{(2i-j) \otimes} \otimes I^{(j-i)\otimes}
\eeqn
with
\beqn{rec}
\begin{array}{lcl}
\mu_{j,i}\nu_i
& = & (a+2-2i) \mu_{j-1,i-1}\nu_{i-1} + (2i-j+1) \mu_{j-1,i}\nu_i \\*[1.5ex]
& = & \big[\mu_{j-1,i-1} + (2i-j+1) \mu_{j-1,i}\big]\nu_i,
\end{array}
\eeqn
where we used the identity
\beqn{nuratio}
\nu_i = (a+2-2i)\nu_{i-1} \tim{ for } i = 1, \ldots, j
\eeqn
to deduce the second equality. Now \req{nablaj-ns-full} gives that
\[
\nabla_\cdot^j \|r\|^a[v]^j
= \bigsum_{i=1}^j \mu_{j,i} \nu_i
\|r\|^{a-j} \, \left(\frac{r^Tv}{\|r\|}\right)^{2i-j} (v^Tv)^{j-i}.
\]
It is then easy to see that the maximum in \req{Tnorm} is achieved for
$v = r/\|r\|$, so that
\beqn{c0}
\| \, \nabla_\cdot^j \| r\|^a \,\|_{[j]}
=\left|\bigsum_{i=1}^j \mu_{j,i}  \nu_i \right| \|r\|^{a-j}
= |\pi_j|\, \|r\|^{a-j}
\eeqn
with
\beqn{pimunu}
\pi_j \eqdef \bigsum_{i=1}^{j}\mu_{j,i}\,\nu_i.
\eeqn
Successively using this definition, \req{rec}, \req{nuratio}
(twice), the identity $\mu_{j-1,j} = 0$ and \req{pimunu} again, we then deduce that
\beqn{pi-ineqs}
\begin{array}{lcl}
\pi_j
& = & \bigsum_{i=1}^{j} \mu_{j-1,i-1}\nu_i + \bigsum_{i=1}^{j} (2i-j+1) \mu_{j-1,i}\nu_i\\*[1.5ex]
& = & \bigsum_{i=1}^{j-1} \mu_{j-1,i}\nu_{i+1} + \bigsum_{i=1}^{j} (2i-j+1) \mu_{j-1,i}\nu_i\\*[1.5ex]
& = & \bigsum_{i=1}^{j-1} \mu_{j-1,i}\big[ \nu_{i+1} + (2i-j+1) \nu_i\big]\\*[1.5ex]
& = & \bigsum_{i=1}^{j-1} \mu_{j-1,i}\big[ (a+2-2(i+1))\nu_i + (2i-j+1) \nu_i\big]\\*[1.5ex]
& = & (a+1-j) \bigsum_{i=1}^{j-1} \mu_{j-1,i}\,\nu_i\\*[1.5ex]
& = & (a+1-j) \pi_{j-1}.
\end{array}
\eeqn
Since $\pi_1 = a$ from the first part of \req{donce}, we
obtain from \req{pi-ineqs} that
\beqn{pij}
\pi_j = \pi(a-j),
\eeqn
which, combined with \req{c0} and \req{pimunu}, gives
\req{dersnorm}. Moreover, \req{pij}, \req{pimunu} and \req{nablaj-ns-full}
give \req{ders} with $\phi_{i,j}= \mu_{j,i}\,\nu_i$. In order to prove
\req{difdersnorm} (where now $\|r\|=1$), we use \req{nablaj-ns-full},
\req{pimunu}, \req{pij} and obtain that
\[
\begin{array}{ll}
\nabla_\cdot^j \|\beta_1r\|^a-\nabla_\cdot^j \|\beta_2r\|^a
&= \bigsum_{i=1}^j \mu_{j,i} \nu_i \,
   \|\beta_1r\|^{a-2i} \, \beta_1^{(2i-j)} r^{(2i-j) \otimes} \otimes
   I^{(j-i)\otimes}\\*[1ex]
& \hspace*{30mm} - \bigsum_{i=1}^j \mu_{j,i} \nu_i \,
   \|\beta_2r\|^{a-2i} \, \beta_2^{(2i-j)} r^{(2i-j) \otimes} \otimes I^{(j-i)\otimes}\\*[4ex]
&= \pi(a-j) \left[\beta_1^{a-j}-\beta_2^{a-j}\right]\,
     \|r\|^{a-2i} \, r^{(2i-j) \otimes} \otimes I^{(j-i)\otimes}\\*[3ex]
&= \pi(a-j) \left[\beta_1^{a-j}-\beta_2^{a-j}\right]\,
      r^{(2i-j) \otimes} \otimes I^{(j-i)\otimes}.
\end{array}
\]
Using \req{Tnorm} again, it is easy to verify that the maximum defining the
norm is achieved for $v=r$ and \req{difdersnorm} then follows from $\|r\|=1$.

\end{document}